\documentclass[a4paper,12pt]{article}
\usepackage{latexsym}
\usepackage{amsmath}
\usepackage{amssymb}
\usepackage{tikz}
\usepackage{color}

\usepackage{theorem}
\newtheorem{theorem}{Theorem}[section]
\newtheorem{proposition}[theorem]{Proposition}
\newtheorem{lemma}[theorem]{Lemma}

\newtheorem{MTA}{Theorem A}
\newtheorem{MTB}{Theorem B}
\newtheorem{MTC}{Theorem C}

\theorembodyfont{\rmfamily}
\newtheorem{proof}{\textmd{\textit{Proof.}}}

\newtheorem{remark}[theorem]{Remark}

\newtheorem{definition}[theorem]{Definition}

\makeatletter

\@addtoreset{equation}{section}
\makeatother

\newcommand{\qedd}{\hfill \Box}

\newcommand{\R}{\ensuremath{\mathbb{R}}}

\newcommand{\Sph}{\ensuremath{\mathbb{S}}}

\def\supp{\mathop{\mathrm{supp}}\nolimits}

\title{ 
Generalized von Mangoldt surfaces of revolution and asymmetric two-spheres of revolution with simple cut locus structure
}
\author{ Minoru Tanaka, Toyohiro Akamatsu, \\ Robert Sinclair and Masaru Yamaguchi }
\date{}
\begin{document}


\maketitle

\begin{abstract}
It is known that if the Gaussian curvature function   along each meridian
on a surface of revolution $(\R^2,dr^2+m(r)^2d\theta^2)$ is decreasing, then
the cut locus of each point of $\theta^{-1}(0)$ is empty or a subarc of the opposite meridian $\theta^{-1}(\pi).$
Such a surface is called a {\it von Mangoldt's surface of revolution}.
A surface of revolution $(\R^2,dr^2+m(r)^2d\theta^2)$ is called a {\it generalized von Mangoldt surface of revolution} if the cut locus of each point of 
$\theta^{-1}(0)$ is empty or a subarc of the opposite meridian $\theta^{-1}(\pi).$

For example, the surface of revolution $(\R^2,dr^2+m_0(r)^2d\theta^2),$ where
$m_0(x)=x/(1+x^2),$ has the same cut locus structure as above and the cut locus of each point in $r^{-1}( (0,\infty) )$ is nonempty.  Note that the Gaussian curvature function is not decreasing along a meridian for this surface.
In this article, we give sufficient conditions for a surface of revolution $(\R^2,dr^2+m(r)^2d\theta^2)$ to be a generalized von Mangoldt surface of revolution.
Moreover, we prove that for any surface of revolution with finite total curvature $c,$ there exists a generalized von Mangoldt surface of revolution with the same total curvature $c$ such that
the Gaussian curvature function  along a meridian is not monotone on $[a,\infty)$ for any $a>0.$
\end{abstract}

\section{Introduction}
As stated in \cite{ASTY}, or \cite{BCST},
determining the structure of the cut locus for a Riemannian manifold
is very difficult, even for a surface of revolution.
In 1994, Hebda \cite{H}  proved that the cut locus of a point $p$ in a complete 2-dimensional Riemannian manifold has a local tree structure and that the cut locus
is locally the image of a Lipschitz map from a subarc  of the unit circle  in the unit tangent plane at $p$
(see also \cite{IT} and \cite[Theorem 4.2.1]{SST}).
In general, one cannot determine the structure of the cut locus in more detail for even a 2-dimensional Riemannian manifold.
 In fact, Gluck and Singer \cite{GS}
constructed a 2-sphere of revolution with positive Gaussian curvature which admits
a non-triangulable cut locus.
Fortunately, the cut locus structure has been determined for very familiar surfaces in Euclidean space, such as
paraboloids, hyperboloids,  ellipsoids and tori of revolution (see \cite{E,GMST, IK1, IK2, ST2}).

These familiar surfaces, which are of revolution, are very useful as model surfaces, for obtaining global structure theorems of
Riemannian manifolds, if the cut locus structure is simple.
For example, special 2-spheres of revolution have been employed as model surfaces to obtain various sphere theorems
(see \cite{B, IMS, K}, for example).

Let $\gamma :[0,t_0]\to M$ denote a minimal geodesic segment on a complete Riemannian manifold $M.$ The endpoint $\gamma(t_0)$ is called a {\it cut point} of $p:=\gamma(0)$ along $\gamma$ if no geodesic extension of $\gamma$ (beyond $\gamma(t_0)$) is minimal. The {\it cut locus} of the point $p$ is defined as the set of cut points along all minimal geodesic segments emanating from $p.$


Let us introduce a new Riemannian metric $g$ on the Euclidean plane $\R^2.$
Choose any smooth function $m:[0,\infty)\to [0,\infty)$ satisfying 
\begin{equation}\label{eq:1.1}
m(x)> 0 \quad \mbox{for \:\:all} \quad x>0.
\end{equation}
Then, the function $m$ can be used to define a Riemannian metric $g$ on ${ \R}^2\setminus\{o\}$
by $g=dr^2+m(r)^2d\theta^2,$
where $(r,\theta)$ denotes polar coordinates about the origin $o$ of the Euclidean plane $\R^2.$
If $m$ is extendable to a smooth odd function around $0,$ 
and satisfies \begin{equation}\label{eq:1.2}
m(0)=0 \quad \mbox{and}\quad m'(0)=1,
\end{equation}
then the Riemannian metric 
$g=dr^2+m(r)^2d\theta^2$
 is extendable to a smooth one on the entirety of $\R^2$
 (see \cite[Theorem 7.1.1]{SST}).

In this article, the Riemannian manifold $(\R^2,g)$ with (smooth) Riemannian metric $g=dr^2+m(r)^2 d\theta^2$  is called a {\it surface  of revolution} if the function $m$ satisfies \eqref{eq:1.1}, \eqref{eq:1.2} and is extendable to a smooth odd function around 0.
The Euclidean plane is a typical example of a surface of revolution with the metric
$dr^2+r^2 d\theta^2.$

A surface of revolution $(\R^2,dr^2+m(r)^2d\theta^2)$ is called a {\it von Mangoldt's surface of revolution} if the Gaussian curvature  is decreasing along each meridian.
Here, each geodesic emanating from the origin $o$ is called a {\it meridian} of the surface of revolution, and the origin $o$ is called the  {\it vertex} of the surface.
Moreover, a curve $r=r_0$ for each positive number $r_0$ is called a {\it parallel.}

Paraboloids of revolution and 2-sheeted hyperboloids of revolution are  typical examples of von Mangoldt's surfaces of revolution (see \cite{M}).
Any von Mangoldt's surface of revolution has a simple cut locus structure, i.e.,
it was proven in \cite{T2} that
the cut locus of a point $p\in \theta^{-1}(0)$ is empty or a subarc of the meridian $\theta^{-1}(\pi)$ which is opposite to the point $p.$

The cut locus of the vertex is always empty for surfaces of revolution, and 
the cut locus of any point distinct from the vertex is nonempty on paraboloids of revolution.
On the other hand, the cut locus of a point $p$ is empty
on the 2-sheeted hyperboloid $z=\sqrt{x^2+y^2+1}$ in Euclidean space 
if the point $p$ is sufficiently close to  the vertex,  i.e., the point $(0,0,1).$
In general, it was proven in \cite[Main Theorem]{T1} that a surface of revolution 
$(\R^2, dr^2+m(r)^2d\theta^2)$ admits uncountably many points whose cut loci are empty if and only if $\liminf_{x\to\infty}m(x)>0,$ and $\int^\infty_11/m(x)^2dx<\infty$
(see \cite[Example 7.3.1, Example 7.3.2]{SST}, \cite[Section III]{M} ).

\begin{definition}\label{def1.1}
A  surface of revolution $(\R^2,dr^2+m(r)^2d\theta^2)$ homeomorphic to the Euclidean plane $\R^2$ is called a {\it generalized von Mangoldt surface of revolution} if the cut locus of each point $p\in\theta^{-1}(0)$ is empty or a subarc of the meridian $\theta^{-1}(\pi).$
\end{definition}

A surface of revolution with nonpositive Gaussian curvature  is a (trivial) generalized von Mangoldt surface of revolution, since the cut locus of any point  is  empty.
Von Mangoldt's surfaces of revolution have served as useful model surfaces.
In fact,
by bounding the sectional curvature of complete open Riemannian manifolds from below by the Gaussian curvature of a von Mangoldt's surface of revolution,
various  global structure theorems have been obtained (see \cite{BCI, KT2} for example.) Hence,  generalized von Mangoldt surfaces can play a very important role as model surfaces
in the global structure theorems mentioned above.


In this article, we will prove  the following three main theorems.


\begin{MTA}
Let $M:=(\R^2,dr^2+m(r)^2d\theta^2)$ be a surface of revolution.
If the function $m$ satisfies the following four properties, then the surface $M$ is a
(non-trivial)  generalized von Mangoldt surface of revolution.
\begin{description}
\item (M.1)\quad\quad\quad 
For some $r_1\in(1,\infty],$   $m'(x)<0$ on $(1,r_1)$ and 
$m'(x)>0$ on $[0,1).$  
\item
\quad\quad\quad\quad\quad \:\;When $r_1<\infty,$   $m'(r_1)=0,$ and $m''(x)\geq 0$ for all $x\geq r_1.$ 

\item
 (M.2)\quad\quad\quad  The function\:\: $-m''(x)/m(x)$ is  decreasing  on $ (0,1) .$
\item
(M.3)\quad\quad\quad  The function $\psi(\nu)$ is decreasing on $( m(r_1),m(1)),$ where
\item\quad\quad\quad \quad\quad\:
$m(r_1):=\lim_{x\to\infty}m(x)$ when $r_1=\infty.$
\item 
\quad \quad \quad \quad \quad \: Here 
$\psi(\nu):=\int^{\eta(\nu)}_{1}
\frac{2\nu}{m(x)\sqrt{m(x)^2-\nu^2}}dx$ and $\eta(\nu):=(m|_{[1,r_1)})^{-1}(\nu).$

\item
(M.4)\quad\quad\quad 
For any $x\in(1,r_1),$
$1+(\xi\circ m)'(x) \geq 0$ holds, where $\xi(\nu):=(m|_{[0,1]})^{-1}(\nu).$

\end{description}
\end{MTA}

\begin{remark} The  function $m(x)$ is strictly decreasing on $[1,\infty)$ when $r_1=\infty.$
Thus,  the limit $\lim_{x\to \infty}m(x)$ exists and $m(1)>\lim_{x\to \infty}m(x)\geq 0,$
since $m(x)>0 $ for all $x>0.$
Hence, 
  $m(1)>m(r_1)\geq 0,$ when $r_1=\infty.$
\end{remark}

The property (M.4) in  Theorem A is not geometrical, but 
we get the following  Theorem B as a corollary to  Theorem A 
under more  geometrical assumptions.
\begin{MTB}
Let $M:=(\R^2,dr^2+m(r)^2d\theta^2)$ be a surface of revolution.
If the function $m$ satisfies the properties (M.5) and (M.6), then the surface $M$ is a (non-trivial)  generalized von Mangoldt surface of revolution.

\begin{description}
 \item (M.5)
\: There exists a positive number $\lambda$ such that the function $m'(x)$ is positive on $[0,\lambda),$
\item \quad \quad \quad \: and $m'(\lambda)=0>m''(\lambda).$
\item (M.6) \: There exists a number   $r_{dc}\in(\lambda,\infty)$ such that  
 the function $-m''(x)/m(x)$ is 
\item \quad \quad \quad \: decreasing on $(0,r_{dc}),$
and   $m''(x)\geq 0$ on $[r_{dc},\infty).$ 
\end{description}
\end{MTB}

\begin{MTC}
For each surface of revolution $\widetilde M=(\R^2,dr^2+\widetilde{m}(r)^2d\theta^2)$ with finite total curvature $c(\widetilde M),$
there exists a generalized von Mangoldt surface of revolution $M$ with the same total curvature as that of $\widetilde M$ such that 
the Gaussian curvature function along each meridian is not monotone on $[a,\infty)$
for any $a>0.$ 
\end{MTC}
\begin{remark}
The finite total curvature does not always imply the boundedness of the Gaussian curvature. In fact, such a surface of revolution
with unbounded Gaussian curvature was constructed in \cite{TK}.
\end{remark}
\begin{remark}
It was conjectured and experimentally illustrated in \cite{ST1} that if a surface of revolution $M$ homeomorphic to the Euclidean plane or a 2-sphere is a union of 
at most $k$ closed domains $D_1,\dots,D_k$ bounded by at most two parallels, such that
the Gaussian curvature function of $M$ along each meridian is monotone in each domain $D_i,$ then the cut locus of each point of $M$ has at most $k$ endpoints.
Thus, Theorem C above and \cite[Theorem A]{ASTY} show that monotonicity of the Gaussian curvature is not a broad enough sufficient condition for a surface of revolution to have a simple cut locus structure.
\end{remark}

In Sections 2, and  3,  Theorem A and   Theorem B will be proven respectively,
in Section 4, some examples of non-trivial generalized von Mangoldt surfaces of revolution are introduced, and in Section 5,   Theorem C will be proven.
In the final Section, as a by-product of some methods used in the proof of 
 Theorem A, a family of {\it asymmetric} 2-spheres of revolution with simple cut locus structure will be introduced. This example is the first one of its kind. In fact, all previously known examples with simple cut locus structure have a reflective symmetry with respect to the equator.




\section{Sufficient conditions}
For basic tools on  surfaces of revolution we refer to \cite{SST}.
Let 
$M:=(\R^2,dr^2+m(r)^2d\theta^2)$ be a surface of revolution such that
 the function $m$ satisfies the above four properties (M.1),...,(M.4).
Note that the function $m$ satisfies \eqref{eq:1.1}, \eqref{eq:1.2} and is extendable to a smooth odd one around 0.
It is not difficult to compute  the Gaussian curvature $G$  of the surface of revolution $M.$ In fact, the curvature $G$ equals 
$$G(q)=\frac{-m''(r(q))}{m(r(q))}$$
at each point $q\in M\setminus\{o\}.$

Therefore, by the property (M.2), the Gaussian curvature function along each meridian is decreasing on the domain 
$r^{-1}( (0,1) ).$
Since $m'(1)=0$ due to the property (M.1), by \cite[Lemma 7.1.4]{SST}, the parallel $r=1$ is a geodesic. Note that (M.1) implies that $m'\geq 0$ on $[r_1,\infty),$ and the Gaussian curvature is nonpositive on $r^{-1}([r_1,\infty)).$
Let $\gamma_{\nu}$ and $\beta_{\nu}$ denote  the unit speed geodesics emanating from a common point $q_0\in\theta^{-1}(0)$ on the parallel  $r=1$ with Clairaut constant $\nu\in(0,m(1)]$
such that 
$\angle((\frac{\partial}{\partial r})_{q_0},\gamma_{\nu}{}'(0))\geq \pi/2$ and $\angle((\frac{\partial}{\partial r})_{q_0},\beta_{\nu}{}'(0))\leq \pi/2$ respectively,
where $\angle(\cdot,\cdot)$ denotes the angle made by  two vectors.

Now let us review the behavior of geodesics $\gamma_\nu$ and $\beta_\nu.$
For each $\nu\in(0,m(1))$ (respectively
 $\nu\in(m(r_1),m(1)),$
the geodesic $\gamma_{\nu}$
(respectively $\beta_\nu$) intersects the parallel $r=1$ again
at a point  $\gamma_{\nu}(a(\nu))$ (respectively $\beta_{\nu}(b(\nu))$),
after the geodesic is  tangent to the parallel $r=\xi(\nu):=(m|_{[0,1]})^{-1}(\nu)$
(respectively  $r=\eta(\nu):=(m|_{[1,r_1)})^{-1}(\nu)$).
Hence by the proof of Proposition 7.1.3 in \cite{SST},
we obtain
that 
$$
\theta(\gamma_{\nu}(a(\nu)))=\varphi(\nu),
\quad \theta(\beta_{\nu}(b(\nu)))=\psi(\nu),  
$$
when $\gamma_\nu([0,a(\nu)])\subset \theta^{-1}([0,2\pi) )$ and  $\beta_\nu([0,a(\nu)])\subset \theta^{-1}([0,2\pi) )$ respectively, 
where 
\begin{equation}\label{eq:2.1}
\varphi(\nu):=\int_{\xi(\nu)}^{1}
\frac{2\nu}{m(x)\sqrt{m(x)^2-\nu^2}}dx\quad\mbox{and}\quad
\psi(\nu):=\int^{\eta(\nu)}_{1}
\frac{2\nu}{m(x)\sqrt{m(x)^2-\nu^2}}dx.
\end{equation}
Since it is possible that $\varphi(\nu)\geq 2\pi$ or  $\psi(\nu)\geq 2\pi$  for some $\nu>0$, the domain of $\theta$ cannot cover the geodesic segment $\gamma_\nu([0,a(\nu)])$ or $\beta_\nu([0,b(\nu)])$  in this case.
For technical reasons, to avoid such inconvenience,
we need to introduce the universal covering space $\widetilde M:=( (0,\infty)\times\R, d\tilde r^2+m(\tilde r)^2d\tilde\theta^2) $ of the space $M\setminus\{o\}.$

For each $\nu\in(0,m(1)),$ let $\tilde \gamma_\nu$ and $\tilde\beta_\nu$
denote the unit speed geodesics emanating from a common point $\tilde q_0$ on the arc $\tilde r=1$ with Clairaut constant $\nu$ such that 
$\angle((\frac{\partial}{\partial \tilde r})_{\tilde q_0},\tilde\gamma_{\nu}{}'(0))\geq \pi/2$ and $\angle((\frac{\partial}{\partial \tilde r})_{\tilde q_0},\tilde \beta_{\nu}{}'(0))\leq \pi/2$ respectively.
From now on we assume that $\pi(\tilde q_0)=q_0,$ 
 where $\pi$ denotes the canonical projection from $\tilde M$ onto $M\setminus\{o\}.$
It is clear that $\pi\circ\tilde\beta_\nu=\beta_\nu,\: \pi\circ\tilde\gamma_\nu=\gamma_\nu,$   and 
$\beta_\nu(t_c),$ (respectively $\gamma_\nu(t_c)),$ where  $t_c>0,$ is conjugate to $ q_0$ along the geodesic $\beta_\nu$ (respectively $\gamma_\nu$) if and only if $\tilde\beta_\nu(t_c)$ (respectively $\tilde\gamma_\nu(t_c))$  is conjugate to $\tilde q_0$ along the geodesic $\tilde\beta_\nu$ (respectively $\tilde\gamma_\nu$). 
Then, we get, by \cite[Proposition 7.1.3]{SST},
$$
\tilde\theta(\tilde\gamma_{\nu}(a(\nu)))=\varphi(\nu)\:\: {\rm and} \:\:
\tilde\theta(\tilde\beta_\nu (b(\nu) ) )=\psi(\nu).
$$

The functions $\varphi(\nu)$ and $\psi(\nu)$ are called the {\it lower half  period function}
and  the {\it upper half period function},
respectively. By (M.1), the parallel $r=r_1$ is a geodesic, 
$\beta_{m(r_1)}$ is not tangent to any parallel, and the points $\beta_{m(r_1)}(s)$ converge to the parallel $r=r_1$ as $s$ goes to infinity when $r_1$ is finite (see \cite [Fig. 7.1.2 and Lemma 7.1.7]{SST}).
For each $\nu\in(0,m(r_1)),$ $\beta_\nu$ is not tangent to any parallel due to the property  (M.1) and by the Clairaut relation (\cite[ (7.1.20)]{SST}). Note that $\lim_{t\to\infty}r(\beta_\nu(t))=\infty$ if $\nu\in(0,m(r_1)).$
\begin{lemma}\label{lem2.1}
The  function $\varphi(\nu)$ is increasing on $(0,m(1)),$ $\varphi(\nu)\geq\pi$ for all $\nu\in(0,m(1)),$ and   $\pi\leq\lim_{\nu\uparrow m(1)}\varphi(\nu)=\lim_{\nu\uparrow m(1)}\psi(\nu).$   
The subarc $r^{-1}(1)\cap\theta^{-1}([0,\pi])$ of the parallel $r=1$ is minimal.
Furthermore, for each $\nu\in(0,m(1)],$ $\tilde\gamma_\nu,$ the  lift of $\gamma_\nu,$
has a conjugate point of $\tilde q_0$ in $\tilde r^{-1}((0,1]),$ but  not  in $\tilde \theta^{-1}( (0,\pi) ).$ 
\end{lemma}
\begin{proof}
From (M.2) and \cite[Lemma 2.2]
{ASTY},
it is clear that $\varphi(\nu)$ is increasing.
Note that
the geometrical proof for this claim is given in
  \cite[Lemma 3.4] {ST2}.
For each $\nu\in(m(r_1),m(1)),$ both  geodesics $\gamma_{\nu}$ and $\beta_{\nu}$ intersect the parallel $r=1$ again at $\gamma_{\nu}(a(\nu))$ and $\beta_{\nu}(b(\nu))$, respectively.  Both intersections $\gamma_{\nu}(a(\nu))$ and $\beta_{\nu}(b(\nu))$
converge to the first conjugate point of $q_0$ along the parallel $r=1$ as $\nu$ tends to $m(1).$  Hence,  by (M.3), $\lim_{\nu\uparrow m(1)}\psi(\nu)= \lim_{\nu\uparrow m(1)}\varphi(\nu) \geq \lim_{\nu\downarrow 0}\varphi(\nu)  =\pi.$ This implies that
the subarc $r^{-1}(1)\cap\theta^{-1}([0,\pi])$ is minimal and that
 the geodesic $\tilde\gamma_{m(1)}=\tilde\beta_{m(1)}$ has no conjugate point of $\tilde q_0$ in $\tilde\theta^{-1}( (0,\pi) ).$  The last statement of our lemma  is clear from the proof of Lemma 3.1 in \cite{ST2}.
$\qedd$
\end{proof}

\begin{lemma} \label{lem2.2}
Suppose that $\nu\in(m(r_1),m(1)),$ where   $m(r_1)$ denotes the number defined in  (M.3).
 Then for any $t\in(0,b(\nu)), $
$\beta_{\nu}(t)$ is not conjugate to $q_0$ along $\beta_{\nu}.$
\end{lemma}
\begin{proof}
Suppose that for some $\nu_0\in(m(r_1),m(1)),$
there exists a conjugate point 
$\beta_{\nu_0}(t_c),$ where $t_c\in(0,b(\nu_0)),$  of $q_0$ along $\beta_{\nu_0}.$
From \cite [Proposition 7.2.2, and Corollary 7.2.1]{SST} the conjugate point 
$\beta_{\nu_0}(t_c)$ appears after $ \beta_{\nu_0}$ is tangent to the parallel 
$r=\eta(\nu_0):=(m|_{[1,r_1)})^{-1}(\nu_0).$
Since $t_c\in(0,b(\nu_0)),$ we get
$r(\beta_{\nu_0}(t_c))>1.$ From \cite[Proposition 7.2.3 and   (7.2.28)]{SST},  it follows that 
\begin{equation*}
\left.\frac{\partial}{\partial\nu}\right|_{\nu_0}
\theta(\beta_{\nu} (t_c))=
\left. \frac{\partial}{\partial\nu}\right|_{\nu_0}
r(\beta_{\nu}(t_c))=0.
\end{equation*}
By applying \cite[(7.1.24)]{SST}
we obtain
$$ \theta(\beta_\nu(t)) -\theta(\beta_\nu(0))
\equiv
\psi(\nu)-\int_1^{r(\beta_\nu(t)) }
\frac{\nu}{m(x)\sqrt{m(x)^2-\nu^2}   }dx  \pmod {2\pi}
$$
for
each $\nu\in(m(r_1),m(1))$ and $t\in(b(\nu)/2,b(\nu)).$
Hence, 
$$0=\psi'(\nu_0)-\int_1^{r(\beta_{\nu_0}(t_c)  )}\frac{m(x)}{\sqrt{m(x)^2-\nu_0^2}^3}dx.$$
This equation contradicts  (M.3), since $r(\beta_{\nu_0}(t_c))>1.$

$\qedd$
\end{proof}
\begin{definition}\label{def2.3}
For each point $q\in r^{-1}( (0,\infty) )\cap\theta^{-1}(0),$ and $\nu\in(0,m( r(q) )],$
let  $\gamma_{\nu}{}^{(q)}$ and $\beta_{\nu}{}^{(q)}$ denote  the unit speed geodesics emanating from the common point $q$ with Clairaut constant $\nu$
such that 
$\angle((\frac{\partial}{\partial r})_{q},(\gamma_{\nu}{}^{(q)}){}'(0))\geq \pi/2$ and $\angle((\frac{\partial}{\partial r})_{q},(\beta_{\nu}{}^{(q)}){}'(0))\leq \pi/2$ respectively.
\end{definition}
\begin{lemma}\label{lem2.4}
For each point $q\in r^{-1}( (0,\infty) )\cap \theta^{-1}(0),$
the geodesic segment
$\alpha^{(q)}|_{[0,s]},$
where
$\alpha^{(q)}:=\gamma_{\nu}^{(q)}=\beta_{\nu}^{(q)}$ for $\nu=m(r(q)),$ has no conjugate point of $q$ if $\alpha^{(q)}([0,s])\subset\theta^{-1}([0,\pi)).$ 
\end{lemma}
\begin{proof}
By Lemma \ref{lem2.1}, we may assume that $r(q)\ne 1.$
Suppose that  $q\in r^{-1}([ r_1,\infty)).$ Then, $ \alpha^{(q)}$ stays in $r^{-1}([r_1,\infty))$ by the Clairaut relation (\cite[Theorem 7.1.2]{SST}) and the property (M.1). Hence, by (M.1), $\alpha^{(q)}$ has no conjugate point of $q.$

Suppose that $q\in r^{-1}( (1,r_1) ).$  In   this case,
 we need to introduce  a lift  $\tilde\alpha^{(q)} $  (on the universal covering space $\widetilde M$) of $\alpha^{(q)},$ which is a geodesic  with Clairaut constant $\nu_q:=m(r(q))$ satisfying  $\pi\circ\tilde\alpha^{(q)}=\alpha^{(q)}.$ 
Then $\tilde \alpha^{(q)}$ 
 is tangent to the arc $\tilde r=\eta(\nu_q)$   at $\tilde q:=\tilde\alpha^{(q)}(0)$ and  to the arc  $\tilde  r=\xi(\nu_q)$ at a point $\tilde q_c.$ The point $\tilde q_c$ is the first conjugate point of $\tilde q$
 along $\tilde\alpha^{(q)}$ by \cite[Lemma 2.8]{ST2} and  \cite[Corollary 7.2.1]{SST}.
Therefore, by Lemma \ref{lem2.1}, 
$\tilde\theta(\tilde q_c)-\tilde\theta(\tilde q)=(\varphi(\nu_q)+\psi(\nu_q))/2\geq\pi,$ and $\alpha^{(q)}|_{[0,s]}$ has no conjugate point of $q$  along $\alpha^{(q)} $ if  $\alpha^{(q)}([0,s]) \subset \theta^{-1}([0,\pi)).$

Suppose that $q\in r^{-1}( (0,1) ).$ If $\nu_q>m(r_1),$ then $\alpha^{(q)}$ is tangent to the parallel $r=\xi(\nu_q)$ at $q$ and to the parallel $r=\eta(\nu_q).$ 
Hence, by the same argument above, we can prove that $\alpha^{(q)}|_{[0,s]}$ has no conjugate point of $q$  along $\alpha^{(q)} $ if  $\alpha^{(q)}([0,s]) \subset \theta^{-1}([0,\pi)).$
If $\nu_q\leq m(r_1),$ then $\alpha^{(q)}$ is not tangent to any parallel except at the point $q.$  By \cite[Corollary 7.2.1]{SST}, there is no conjugate point of $q$ along $\alpha^{(q)}.$
$\qedd$
\end{proof}

\begin{lemma}\label{lem2.5}
Suppose that    $q\in r^{-1}( (0,r_1) ).$ If $m(r_1)>0,$
then for each $\nu\in(0,m(r_1)  ]\cap(0,m(r(q))),$
the geodesic  $\beta_\nu{}^{(q)}$ 
has no conjugate point of 
$q.$ 
\end{lemma}
\begin{proof}
Since  $\nu\leq m(r_1),$ 
 $\beta_\nu^{(q)}$ is not tangent to any parallel, and hence the geodesic has no conjugate point of $q.$ 
Note that   
 $\lim_{t\to\infty}r( \beta_{m(r_1)}^{(q)}(t) )=r_1$ if $r_1<\infty$  and $\lim_{t\to\infty}r(\beta_\nu^{(q)}(t))=\infty$, if $\nu<m(r_1).$

$\qedd$
\end{proof}

\begin{lemma}\label{lem2.6}
Suppose that $q\in r^{-1}([r_1,\infty)).$
Then,
for each $\nu\in(0,m(r(q))),$ the geodesic  $\beta_\nu^{(q)}$ has no conjugate point  of
$q.$
\end{lemma}
\begin{proof}
By the property (M.1), $m\geq m(r_1)$ on $[r_1,\infty).$
Since the geodesic $\beta_\nu^{(q)}(t)$  is not tangent to any parallel by the Clairaut relation (see \cite[Theorem 7.2.2]{SST}),
the geodesic  has no conjugate point of $q$ by \cite[Proposition 7.2.1]{SST}.

$\qedd$
\end{proof}
\begin{lemma}\label{lem2.7}
Suppose that $q\in r^{-1}([r_1,\infty)).$
Then, for each $\nu\in [m(r_1),m(r(q) ) ],$
the geodesic  $\gamma_\nu^{(q)}$ has no conjugate point of $q.$
\begin{proof}
If $r(q)>r_1,$ then
$\gamma_\nu^{(q)}$ does not intersect the parallel $r=r_1$ by the Clairaut relation (\cite[(7.1.20)]{SST}).
Thus, 
it is clear from (M.1) that $\gamma_\nu^{(q)}$ has no conjugate point of $q.$
If $r(q)=r_1,$ then $m(r(q))=m(r_1)$ and $\nu=m(r(q))=m(r_1)$, since $\nu\in[m(r_1),m(r(q))].$ This implies that the  geodesic $\gamma_\nu^{(q)}$ is the parallel $r=r_1.$
By (M.1), the geodesic has no conjugate point of $q.$
$\qedd$
\end{proof}

\end{lemma}

We will prove the following two propositions which are crucial for the proof of 
{ Theorem A.}

\begin{proposition}\label{prop2.8}
 For each $q\in r^{-1}([1,r_1))\cap \theta^{-1}(0)$
and $\nu\in(0,m(r(q) )  ],$
the geodesic segment 
$\gamma_{\nu}^{(q)}|_{[0,s]}$ has no conjugate point of $q$ if
$\gamma_\nu^{(q)}([0,s])\subset\theta^{-1}([0,\pi)).$
\end{proposition}

\begin{proposition}\label{prop2.9}
For each $q\in r^{-1}([r_1,\infty))\cap\theta^{-1}(0),$ and $\nu\in(0,m(r(q)) ], $
the geodesic segment $\gamma_\nu^{(q)}|_{[0,s]}$ has no conjugate point of $q$ if $\gamma_\nu^{(q)}([0,s])\subset\theta^{-1}([0,\pi)).$
\end{proposition}


We need the following series of  lemmas (Lemmas \ref{lem2.10},..., \ref{lem2.17})
for proving Propositions \ref{prop2.8} and  \ref{prop2.9}. Then, 
one can prove  Proposition \ref{prop2.8} (respectively Proposition \ref{prop2.9})
by combining Lemmas \ref{lem2.4}, \ref{lem2.11}, \ref{lem2.15}, and \ref{lem2.17} (respectively Lemmas \ref{lem2.7} and \ref{lem2.15}).

\begin{lemma}\label{lem2.10}
For each $\nu_0\in(m(r_1),m(1)),$ and each point $q\in r^{-1}([1,\eta(\nu_0))),$ 
$\gamma_{\nu_0}^{(q)}$ has a conjugate point of $q$ and 
\begin{equation}\label{eq:2.2}
\varphi'(\nu_0)+\int_1^{r(q)}f_\nu(x,\nu_0)dx-\int_{r(q_c)}^1f_\nu(x,\nu_0)dx=0
\end{equation}
holds, where 
$$ 
f(x,\nu):=
\nu\left(m(x)\sqrt{m(x)^2-\nu^2}\right)^{-1},\:
f_\nu(x,\nu):=\frac{\partial f}{\partial \nu} (x,\nu)=m(x)\left({m(x)^2-\nu^2} \right)^{-3/2}, $$ and $q_c$ denotes the first conjugate point of $q$ along $\gamma_{\nu_0}^{(q)}.$
\end{lemma}

\begin{proof}
Choose any $\nu_0\in(m(r_1),m(1)),$ and any point $q\in r^{-1}([1,\eta(\nu_0))),$ and fix them.
Write
$\alpha_\nu:=\gamma_\nu^{(q)}$ for simplicity for each $\nu\in( m(r_1),m(r(q)) )$ in this proof.
The geodesic $\alpha_{\nu_0}$ intersects the parallel
$r=1$ at a point $\alpha_{\nu_0}(t(\nu_0)),$  which denotes the first intersection, and is tangent to the parallel $r=\xi(\nu_0)$.

Since the subarc $\gamma_{\nu_0}$ of the geodesic $\alpha_{\nu_0}$ emanating from $\alpha_{\nu_0}(t(\nu_0))$ has a conjugate point of $\alpha_{\nu_0}(t(\nu_0))$ in $r^{-1}((0,1])$ by Lemma \ref{lem2.1},
$\alpha_{\nu_0}$ also has a  conjugate point $q_c:=\alpha_{\nu_0}(t_c(\nu_0))$ of $q$ in $ r^{-1}( (0,1] ).$
From now on, we assume that the conjugate point $q_c$ is the first one and write $t_c:=t_c(\nu_0)$ for simplicity.

Since $\alpha_{\nu_0}(t_c)$ is the first conjugate point of $q$ along $\alpha_{\nu_0},$
we get, by \cite[Proposition 7.2.3 and (7.2.28) ]{SST}
\begin{equation}\label{eq:2.3}
\left.\frac{\partial}{\partial\nu}\right|_{\nu_0}\theta(\alpha_\nu(t_c))=\left.\frac{\partial}{\partial\nu}\right|_{\nu_0} r(\alpha_\nu(t_c))=0.
\end{equation}
By applying \cite[(7.1.24)]{SST}, we obtain
\begin{equation}\label{eq:2.4}
\theta(\alpha_\nu(t_c))-\theta(q)\equiv\varphi(\nu)+\int_1^{r(q)}f(x,\nu)dx-\int_{r(\alpha_\nu(t_c))}^1f(x,\nu)dx \pmod{2\pi}
\end{equation}
for any $\nu$ sufficiently close to $\nu_0.$
By \eqref{eq:2.3} and \eqref{eq:2.4},
we get  \eqref{eq:2.2}.
$\qedd$
\end{proof}

Lemma \ref{lem2.10} implies that
for each $t\in[1,\eta(\nu_0)),$ where the number $\nu_0\in(m(r_1),m(1))$
is assumed to be fixed, there exists a number $h(t)\in(\xi(\nu_0),1]$
such that
\begin{equation}\label{eq:2.5}
\varphi'(\nu_0)+\int_{1}^{t}f_{\nu}(x,\nu_0)dx-\int_{h(t)}^{1} f_\nu(x,\nu_0)dx=0,
\end{equation}
and, for each point $q$ on the parallel $r=t\in[1,\eta(\nu_0)),$ there exists a first conjugate point $q_c$ on the parallel $r=h(t)$ of $q$ along $\alpha_{\nu_0}.$

By \eqref{eq:2.5}, $h(t)$ is differentiable on $(1,\eta(\nu_0)),$ and
\begin{equation}\label{eq:2.6}
h'(t)=-{f_\nu(t,\nu_0)}/{f_\nu(h(t),\nu_0)}.
\end{equation}
From the proof of \cite[Proposition 7.1.3]{SST}, it follows that 
\begin{equation}\label{eq:2.7}
\theta(\alpha_{\nu_0}(t_c))-\theta(q)=\Phi_{\nu_0}(r(q),r(\alpha_{\nu_0}(t_c)))=\Phi_{\nu_0}(r(q),h(r(q)))
\end{equation}
holds when $\alpha_{\nu_0}( [0,t_c])\subset \theta^{-1}([0,2\pi))$, where
\begin{equation}\label{eq:2.8}
\Phi_{\nu_0}(u,v):=\varphi(\nu_0)+\int_1^{u}f(x,\nu_0)dx-\int^1_{v}f(x,\nu_0)dx
\end{equation}
for $u\in[1,\eta(\nu_0)),$ and $v\in(\xi(\nu_0),1].$
If we introduce a lift $\tilde\alpha_{\nu_0}$ of $\alpha_{\nu_0}$ on the universal covering space $\widetilde M,$
then $\tilde\alpha_{\nu_0}(t_c)$ is a conjugate point of $\tilde q:=\tilde\alpha_{\nu_0}(0) $ along $\tilde\alpha_{\nu_0},$
and 
 \eqref{eq:2.7} can be expressed  as 
$$   
\tilde\theta(\tilde\alpha_{\nu_0}(t_c))-\tilde\theta(\tilde q)=\Phi_{\nu_0}(r(q),r(\alpha_{\nu_0}(t_c)))=\Phi_{\nu_0}(r(q),h(r(q))).$$

\begin{lemma}\label{lem2.11}
Suppose that  $\nu_0\in(m(r_1),m(1)).$ Then,
for all $t\in[1,\eta(\nu_0)),$
$$\Phi_{\nu_0}(t,h(t))\geq\pi$$ holds. Therefore, for each point $q\in r^{-1}([1,r_1))$ and $\nu_0\in(m(r_1),m(r(q))),$ the geodesic segment $\gamma_{\nu_0}^{(q)}|_{[0,s]}$ has no conjugate point of $q$ if $\gamma_{\nu_0}^{(q)}([0,s])\subset\theta^{-1}([0,\pi)).$
\end{lemma}
\begin{proof}
Suppose that $\Phi_{\nu_0}(t,h(t))<\pi$ for some $t\in(1,\eta(\nu_0)).$ Then
$\Phi_{\nu_0}(t,h(t))$ attains a minimum ($<\pi$) at $t_0\in(1,\eta(\nu_0)),$
since $\lim_{t\to\eta(\nu_0)}\Phi_{\nu_0}(t,h(t))=(\varphi(\nu_0)+\psi(\nu_0))/2\geq\pi$ by   Lemma \ref{lem2.1}, and since 
$\Phi_{\nu_0}(1,h(1))=\tilde\theta(\tilde\gamma_{\nu_0}(t_c(\nu_0)))-\tilde\theta(\tilde\gamma_{\nu_0}(0))\geq\pi$ by Lemma \ref{lem2.1}, where $\tilde\gamma_{\nu_0}$ denotes the  lift of $\gamma_{\nu_0}$ in Lemma \ref{lem2.1} and 
$\tilde\gamma_{\nu_0}(t_c(\nu_0))$ denotes the first conjugate point of $\tilde\gamma_{\nu_0}(0)$ along the geodesic $\tilde\gamma_{\nu_0}.$ 
By \eqref{eq:2.6} and \eqref{eq:2.8}, 
we obtain
$$\frac{d\Phi_{\nu_0}(t,h(t))}{ dt}=f_{\nu}(t,{\nu_0})\left(( f/f_{\nu})(t,\nu_0)-  (f/f_{\nu})(h(t),\nu_0)\right).$$
Since $(f/f_{\nu})(t,\nu)=\nu-\nu^3/m(t)^2,$
we get
\begin{equation*}\label{eq:2.9}
\left.\frac{d}{dt}\right|_{t_0}\Phi_{\nu_0}(t,h(t))=\nu_0^3f_{\nu}(t_0,\nu_0) (m(h(t_0))^{-2}-m(t_0)^{-2}).
\end{equation*}
Since $\left.\frac{d}{dt}\right|_{t=t_0}\Phi_{\nu_0}(t,h(t))=0,$
we obtain $m(t_0)=m(h(t_0)),$ and $\pi>\Phi_{\nu_0}(t_0,h(t_0))=\Phi_{\nu_0}(t_0,\xi\circ m(t_0) ),$
which is a contradiction by the following lemma.
$\qedd$
\end{proof}

\begin{lemma}\label{lem2.12}
For each 
 $\nu_0\in(m(r_1),m(1))$ 

$$\Phi_{\nu_0}(t,\xi\circ m(t))\geq\pi\quad {\rm and} \quad \Psi_{\nu_0}(t,\xi\circ m (t) )\geq\pi$$ hold on  $[1,\eta(\nu_0))$.
Here, 
$$\Psi_{\nu_0}(u,v):=\int^{\eta(\nu_0)}_uf(x,\nu_0)dx+\int_{v}^{\eta(\nu_0)}f(x,\nu_0)dx$$ 
for  $u,v\in(\xi(\nu_0),\eta(\nu_0)).$
\end{lemma}

\begin{proof}
Since $m(t)=m(\xi\circ m(t) )$ for all $t\in[1,\eta(\nu_0)),$
it is easy to check that
$f(t,\nu_0)\cdot(1+(\xi\circ m)'(t))=\frac{d}{dt}\Phi_{\nu_0}(t,\xi\circ m(t))=-\frac{d}{dt}\Psi_{\nu_0}(t,\xi\circ m(t) ).$
Hence, by (M.4), the functions $\Phi_{\nu_0}(t,\xi\circ m(t))$  and $-\Psi_{\nu_0}(t,\xi\circ m(t) )$ are   increasing on $(1,\eta(\nu_0)),$
$\Phi_{\nu_0}(t,\xi\circ m(t))\geq \Phi_{\nu_0}(t,\xi\circ m(t))|_{t=1}=\varphi(\nu_0),$ and $\Psi_{\nu_0}(t,\xi\circ m(t))\geq \lim_{t\uparrow\eta(\nu_0)}\Psi_{\nu_0}(t,\xi\circ m(t))=(\varphi(\nu_0)+\psi(\nu_0))/2$
for all $t\in(1,\eta(\nu_0)).$
From Lemma \ref{lem2.1}, $\Phi_{\nu_0}(t,\xi\circ m(t))\geq \pi$ and $\Psi_{\nu_0}(t,\xi\circ m(t))\geq \pi$ for all $t\in(1,\eta(\nu_0)).$

$\qedd$
\end{proof}

The function $\Phi_{\nu_0}(u,v)$ is well-defined for $\nu_0=m(r_1)$ by
\eqref{eq:2.8} on $[1,r_1)\times(\xi\circ m(r_1),1],$  where $r_1\in(1,\infty]$ is assumed to  satisfy  $m(r_1)>0,$ whereas $\Psi_{\nu_0}(u,v)$ is not well-defined for $\nu_0=m(r_1).$
The proof of the following lemma is the same as that of Lemma \ref{lem2.12}.

\begin{lemma}\label{lem2.13}
Suppose that   $m(r_1)>0.$
Then,
$$\Phi_{m(r_1)}(t,\xi\circ m(t))\geq \pi$$
holds on $[1,r_1).$
\end{lemma}

\begin{lemma}\label{lem2.14}
Suppose that  $m(r_1)>0.$
Then
for  each pair of a point $q\in r^{-1}([1,r_1))$ and a number $\nu_0\in(0,m(r_1)],$ or each pair of a point $q\in r^{-1}([r_1,\infty))$ and a number $\nu_0\in(0,m(r_1)),$
$\gamma_{\nu_0}^{(q)}$ has a conjugate point of $q$ and 
\begin{equation}\label{eq:2.10}
\varphi'(\nu_0)+\int_1^{r(q)}f_\nu(x,\nu_0)dx-\int_{r(q_c)}^1f_\nu(x,\nu_0)dx=0
\end{equation}
holds,  where $q_c$ denotes the first conjugate point of $q$ along $\gamma_{\nu_0}^{(q)}.$
\end{lemma}
\begin{proof}
Choose any pair of a point $q\in r^{-1}([1,r_1))$ and a number $\nu_0\in(0,m(r_1)],$ or any pair of a point $q\in r^{-1}([r_1,\infty))$ and a number $\nu_0\in(0,m(r_1)).$
Then, $\nu_0<m(r(q))$ holds. Thus,
by the Clairaut relation, $\alpha_{\nu_0}:=\gamma_{\nu_0}^{(q)}$ is not tangent to the parallel $r=r(q)$ at $q=\alpha_{\nu_0}(0)$, but  tangent to the parallel 
 $r=\xi(\nu_0).$
Hence, the proof of Lemma \ref{lem2.10} is still valid and we obtain \eqref{eq:2.10}.
Note that $\alpha_{\nu_0}$ is not tangent to any parallel after it is tangent to the parallel $r=\xi(\nu_0)$ in this case.

$\qedd$
\end{proof}

If the number $\nu_0 $ is less than $m(r_1)$ (respectively $\nu_0=m(r_1)$), then
Lemma \ref{lem2.14} implies that for each $t\in[1,\infty)$ (respectively $t\in[1,r_1)$) there exists a number $h(t)\in(\xi(\nu_0),1]$ such that
\begin{equation}\label{eq:2.11}
\varphi'(\nu_0)+\int_{1}^{t}f_{\nu}(x,\nu_0)dx-\int_{h(t)}^{1} f_\nu(x,\nu_0)dx=0,
\end{equation}
and, for each point $q$ on the parallel $r=t\in[1,\infty),$ (respectively $r=t\in[1,r_1))$ there exists a first conjugate point $q_c\in r^{-1}(h(t))$ of $q$ along $\alpha_{\nu_0}.$ 

Furthermore,
by \eqref{eq:2.11}, $h(t)$ is differentiable on $(1,\infty)$ (respectively on $(1,r_1)$), if $\nu_0<m(r_1)$ (respectively  if $\nu_0=m(r_1)$) and
\begin{equation}\label{eq:2.12}
h'(t)=-{f_\nu(t,\nu_0)}/{f_\nu(h(t),\nu_0)}.
\end{equation}
From the proof of  \cite[(7.1.24)]{SST} it follows that 
\begin{equation*}\label{eq:2.13}
\theta(q_c)-\theta(q)=\Phi_{\nu_0}(r(q),r(\alpha_{\nu_0}(t_c)))=\Phi_{\nu_0}(r(q),h(r(q)))
\end{equation*}
holds when $\alpha_{\nu_0}([0,t_c])\subset \theta^{-1}([0,2\pi)),$ where $t_c$ denotes the parameter value satisfying $\alpha_{\nu_0}(t_c)=q_c,$
and 
$$\tilde \theta(\tilde\alpha_{\nu_0}(t_c))-\tilde\theta(\tilde\alpha_{\nu_0}(0))=\Phi_{\nu_0}(r(q),r(\alpha_{\nu_0}(t_c)))=\Phi_{\nu_0}(r(q),h(r(q)))$$
holds,
 where
\begin{equation}\label{eq:2.14}
\Phi_{\nu_0}(u,v):=\varphi(\nu_0)+\int_1^{u}f(x,\nu_0)dx-\int^1_{v}f(x,\nu_0)dx
\end{equation}
for $u\in[1,\infty),$ $v\in(\xi(\nu_0),1]$ if $\nu_0<m(r_1)$ (respectively $u\in[1,r_1),$ $v\in(\xi\circ m(r_1),1]$ if $\nu_0=m(r_1)$), and $\tilde\alpha_{\nu_0}$ denotes a lift of $\alpha_{\nu_0}$ on the universal covering space $\widetilde M.$

\begin{lemma}\label{lem2.15}
Suppose that $\nu_0\in(0,m(r_1) ).$
Then for all $t\in[1,\infty),$
$$\Phi_{\nu_0}(t,h(t))\geq\pi$$ holds.
Therefore, for each point $q\in r^{-1}([1,\infty))$ and $\nu_0\in(0,m(r_1)),$ the geodesic segment $\gamma_{\nu_0}^{(q)}|_{[0,s]}$ has no conjugate point of $q$ if $\gamma_{\nu_0}^{(q)}([0,s])\subset\theta^{-1}([0,\pi)).$

\end{lemma}
\begin{proof}
By \eqref{eq:2.12} and \eqref{eq:2.14},
we get
\begin{equation}\label{eq:2.15}
\frac{d}{dt}\Phi_{\nu_0}(t,h(t))=\nu_0^3f_\nu(t,\nu_0)(m(h(t))^{-2}-m(t)^{-2}).
\end{equation}
Suppose that $\Phi_{\nu_0}(t,h(t))<\pi$ for some $t\in(1,\infty).$ 
We will prove that $\Phi_{\nu_0}(t,h(t))$  attains a minimum at some $t_0\in(1,t_1),$
where
$$t_1:=\sup\{x>1; m(x)<m(1)\}.$$
If $t_1$ is finite, then the function $\frac{d}{dt}\Phi_{\nu_0}(t,h(t))$ is strictly positive on $[t_1,\infty)$ by \eqref{eq:2.15}. Hence, the function attains a minimum at some  $t_0\in(1,t_1)$ when $t_1$ is finite.

Suppose next that  $t_1$  is infinite. Thus, $m(x)\leq m(1)$ for all $x\geq 1,$
and $\int^\infty_1 {1}/{m(x)^2}dx=\infty.$ This implies that
$\lim_{t\to\infty}\Phi_{\nu_0}(t,h(t))=\infty, $ since $\Phi_{\nu_0}(t,h(t))\geq\int_1^t f(x,\nu_0)dx>\int^t_1 {\nu_0}/{m(x)^2}dx.$
Therefore, 
$\Phi_{\nu_0}(t,h(t))$ attains a minimum ($<\pi$) at $t_0\in(1,t_1).$ 
From \eqref{eq:2.15}, we obtain $m(t_0)=m(h(t_0)),$ $h(t_0)=\xi\circ m(t_0)$ and 
$\pi>\Phi_{\nu_0}(t_0,h(t_0) )=\Phi_{\nu_0}(t_0,\xi\circ m(t_0) ),$
which is a contradiction by the following lemma.
$\qedd$
\end{proof}

\begin{lemma}\label{lem2.16}
Suppose that   $\nu_0\in(0,m(r_1)).$ 
Then,
for all $t\in[1,t_1),$ where $t_1=\sup\{x>1; m(x)<m(1)\},$
$$\Phi_{\nu_0}(t,\xi\circ m(t) )\geq\pi$$ holds.
\end{lemma}

\begin{proof}

Since $m(t)=m(\xi\circ m(t) )$ for all $t\in(1,t_1),$
it is easy to check that
\begin{equation}\label{eq:2.16}
\frac{d}{dt}\Phi_{\nu_0}(t,\xi\circ m(t) )=f(t,\nu_0)(1+(\xi\circ m)'(t)).
\end{equation}
Here, we will prove that 
$1+(\xi\circ m)'(x)\geq 0$ on $(1,t_1).$
Then, it is clear from \eqref{eq:2.16} that $\Phi_{\nu_0}(t,\xi\circ m(t) )$ is increasing on $(1,t_1),$ and hence, for all $t\in[1,t_1),$
 $\Phi_{\nu_0}(t,\xi\circ m(t))\geq \Phi_{\nu_0}(1,\xi\circ m(1))=\varphi(\nu_0)\geq \pi$
  by Lemma \ref{lem2.1}. 

By (M.4),  $1+(\xi\circ m)'(x)\geq 0$ on $(1,r_1).$
In particular, $1+(\xi\circ m)'(r_1)\geq 0.$ 
It is sufficient to prove that $(\xi\circ m)'(x)$ is increasing on $[r_1,t_1].$

Since $m'(x)\geq 0$  and $m''(x)\geq 0$ on $(r_1,t_1),$ and since $\xi'(\nu)>0$ and $\xi''(\nu)>0$  on $(0,m(1) ),$ $(\xi\circ m)''(x)=\xi''(m(x))\cdot m'(x)^2+\xi'(m(x))\cdot m''(x)\geq 0$ on $(r_1,t_1).$ 
Hence $(\xi\circ m)'(x)$ is increasing on $(r_1,t_1).$

$\qedd$
\end{proof}
\begin{lemma}\label{lem2.17}
Suppose that $\nu_0=m(r_1)>0.$
Then, for all $t\in[1,r_1),$
$$\Phi_{\nu_0}(t,h(t))\geq\pi$$ holds.
Therefore, for each point $q\in r^{-1}([1,r_1)),$  the geodesic segment $\gamma_{m(r_1)}^{(q)}|_{[0,s]}$ has no conjugate point of $q$ if $\gamma_{m(r_1)}^{(q)}([0,s])\subset\theta^{-1}([0,\pi)).$
\end{lemma}
\begin{proof}
First, we will prove that $\lim_{t\to r_1}\Phi_{\nu_0}(t,h(t))=\infty.$
Suppose that $r_1=\infty.$  Then $m'(x)<0$ on $(1,\infty)$ by (M.1).
This implies that $\lim_{t\to\infty}\Phi_{\nu_0}(t,h(t))=\infty,$ since $\Phi_{\nu_0}(t,h(t)) \geq \int_1^t f(x,m(r_1))dx \geq m(r_1)\int_1^t {1}/{m(x)^2}dx.$

Next, suppose that $r_1$ is finite. By (M.1), $m'(r_1)=0.$
Since $\lim_{x\to r_1}(m(x)-m(r_1))/(x-r_1)^2=m''(r_1)/2$ by l'H\^opital's rule, it is clear that
$ \lim_{t\uparrow r_1}\int_1^tf(x,\nu_0) dx=\infty.$

This property is also geometrically proven by noting that
 the geodesic $\beta_{m(r_1)}^{(q)}(t)$ is not tangent to any parallel  and converges to the parallel $r=r_1$ as $t$ goes to $\infty$ (see \cite[ (7.1.19)]{SST}).

Suppose that $\Phi_{\nu_0}(t,h(t))$ is less than $\pi$ for some $t\in(1,r_1).$
Then, 
$\Phi_{\nu_0}(t,h(t))$ attains a  minimum ($<\pi$) at some number $t=t_0\in(1,r_1),$
since $\Psi_{\nu_0}(1,h(1))\geq\pi$ by Lemma \ref{lem2.1} and since $\lim_{t\to r_1}\Phi_{\nu_0}(t,h(t))=\infty.$
By \eqref{eq:2.15},  we have $h(t_0)=\xi\circ m(t_0). $
Therefore, $\pi>\Phi_{\nu_0}(t_0,h(t_0))=\Phi_{\nu_0}(t_0,\xi\circ m(t_0))\geq\pi$ from Lemma \ref{lem2.13}. This is a contradiction.


$\qedd$
\end{proof}

The proof of 
the following proposition is similar to   that of Lemma \ref{lem2.11}.
\begin{proposition}\label{prop2.18}
For each $\nu_0\in(m(r_1),m(1)),$ and each point $q\in r^{-1}( (\xi(\nu_0),\eta(\nu_0)) )\subset r^{-1}( (0,r_1) ),$
the geodesic segment
$\beta_{\nu_0}^{(q)}|_{[0,s]}$ has no conjugate point of $q$ if  $\beta_{\nu_0}^{(q)}([0,s])\subset \theta^{-1}([0,\pi) ).$
\end{proposition}

\begin{proof}
Choose any $\nu_0\in (m(r_1),m(1))$ and any point $q\in r^{-1}((\xi(\nu_0),\eta(\nu_0))).$  
The geodesic $\beta_{\nu_0}^{(q)}$ is tangent to the parallel $r=\eta(\nu_0),$ and, after its intersection, the geodesic is tangent to the parallel $r=\xi(\nu_0).$
Since each pair of the intersection of the geodesic $\beta_{\nu_0}^{(q)}$ and either the parallel $r=\eta(\nu_0)$ or $r=\xi(\nu_0)$ is a   conjugate pair along $\beta_{\nu_0}^{(q)},$ 
the geodesic has a first conjugate point $q_c$ of $q$ (see \cite[Lemma 2.8]{ST2}).
Therefore, we obtain
\begin{equation}\label{eq:2.17}
\psi'(\nu_0)-\int^{r(q)}_1f_\nu(x,\nu_0)dx+\int_{r(q_c)}^1f_\nu(x,\nu_0)dx=0,
\end{equation}
where $f_\nu(x,\nu)=m(x)(m(x)^2-\nu^2) ^{-3/2}.$

Thus, by \eqref{eq:2.17}, for each $t\in(\xi(\nu_0),\eta(\nu_0)),$
there exists a number $h(t)\in(\xi(\nu_0),\eta(\nu_0) )$ satisfying
\begin{equation}\label{eq:2.18}
\psi'(\nu_0)-\int^{t}_1f_\nu(x,\nu_0)dx+\int_{h(t)}^1f_\nu(x,\nu_0)dx=0,
\end{equation}
and for each point $q$ on the parallel $r=t\in(\xi(\nu_0),\eta(\nu_0)),$ there exists a first conjugate point $q_c$ on the parallel $r=h(t)$ of $q$ along $\beta_{\nu_0}^{(q)}.$
From the inverse function theorem
it follows that $h(t)$ is differentiable on $(\xi(\nu_0),\eta(\nu_0)),$ and 

\begin{equation}\label{eq:2.19}
h'(t)=-{f_\nu(t,\nu_0)}/{f_\nu(h(t),\nu_0)}.
\end{equation} 
From the proof of \cite[ (7.1.24)]{SST} it follows that 
\begin{equation*}\label{eq:2.20}
\theta(q_c)-\theta(q)=\Psi_{\nu_0}(r(q),h(r(q)) )
\end{equation*}
holds when $\beta_{\nu_0}^{(q)}( [0,t_c])\subset \theta^{-1}( [0,2\pi) ),$
where
$\Psi_{\nu_0}$ denotes the function defined in Lemma \ref{lem2.12}
and $\beta_{\nu_0}^{(q)}(t_c)=q_c.$
Note that 
\begin{equation}\label{eq:2.21}
\Psi_{\nu_0}(t,h(t))=\psi(\nu_0)-\int^t_1 f(x,\nu_0)dx+\int^1_{h(t)}f(x,\nu_0)dx
\end{equation}
holds for all $t\in (\xi(\nu_0),\eta(\nu_0)).$ 

Suppose that $\Psi_{\nu_0}(t,h(t)))<\pi$ for some $t\in(\xi(\nu_0),\eta(\nu_0)).$
Since
$\pi\leq(\varphi(\nu_0)+\psi(\nu_0))/2 =\lim_{t\uparrow\eta(\nu_0)}\Psi_{\nu_0}(t,h(t))=\lim_{t\downarrow\xi(\nu_0)}\Psi_{\nu_0}(t,h(t))$ by Lemma \ref{lem2.1}, the function 
$\Psi_{\nu_0}(t,h(t))$ attains a minimum at $t=t_0\in(\xi(\nu_0),\eta(\nu_0))$  and
 $\Psi_{\nu_0}(t_0,h(t_0))<\pi.$
By  Lemmas \ref{lem2.1}, \ref {lem2.2},  \eqref{eq:2.18} and \eqref{eq:2.21},
either
$t_0>1>h(t_0)$  or $t_0<1<h(t_0).$ 
Without loss of generality, we may assume that
$t_0>1>h(t_0).$ 
Since $\left.\frac{d}{dt}\right|_{t_0}\Psi_{\nu_0}(t,h(t))=0,$ we get
$m(t_0)=m(h(t_0))$ by \eqref{eq:2.19} and \eqref{eq:2.21}.  Then, $\xi\circ m(t_0)=h(t_0),$ and $\Psi_{\nu_0}(t_0,h(t_0))=
\Psi_{\nu_0}(t_0,\xi\circ m(t_0))<\pi.$ This is a contradiction by Lemma \ref{lem2.12}.
Therefore, 
the geodesic segment
$\beta_{\nu_0}^{(q)}|_{[0,s]}$ has no conjugate point of $q$ if  $\beta_{\nu_0}^{(q)}([0,s])\subset \theta^{-1}([0,\pi) ).$
$\qedd$
\end{proof}


Here, we will give the proof of  {Theorem A}  stated in the introduction.\\

\noindent
{Proof of Theorem A.}\vspace{3mm}
\\
Suppose that there exists  a point $q\in\theta^{-1}(0)$ having a cut point $q_c\notin\theta^{-1}(\pi).$
We may assume that $q_c\in \theta^{-1}( (0,\pi) )$ and that $q_c$ is a conjugate point of $q$
along a minimal geodesic segment $\alpha$ joining $q$ to $q_c.$
Otherwise, there exists another  minimal geodesic segment $\beta$ joining $q$ to $q_c.$
Then, by imitating the same method used in the proof in \cite[Lemma 3.1]{ST2}, we can find a cut point $q_0$ of $q$ in the domain bounded by $\alpha$ and $\beta$ which is a conjugate point of $q$ along any minimal geodesic segment joining $q$ to $q_0.$
Therefore, there exists a cut point $q_c\in\theta^{-1}( (0,\pi) )$ of $q$ which is  conjugate along a minimal geodesic segment $\alpha$ joining $q$ to $q_c.$
From Lemmas \ref{lem2.4}, \ref{lem2.5}, \ref{lem2.6}, and Proposition \ref{prop2.18} it follows that
$\alpha\ne\beta_{\nu}^{(q)}$ for any $\nu\in(0,m(r(q))].$
Therefore, $\alpha=\gamma_{\nu_0}^{(q)}$ for some $\nu_0\in(0,m(r(q))).$
By Propositions \ref{prop2.8} and \ref{prop2.9},
we obtain $r(q)\in(0,1).$ To be conjugate is a symmetric property, i.e., if $q_c$ is conjugate to $q$ along $\alpha,$ then $q$ is conjugate to $q_c$ along $\alpha.$
Thus, we get $r(q), r(q_c)\in(0,1).$ 


It follows from the proof of Lemma 3.1 in \cite{ST2} that 
$\alpha=\gamma_{\nu_0}^{(q)}$ has no conjugate point of  the point $q\in r^{-1}( (0,1) )$ in $\theta^{-1}( (0,\pi) ).$ This is a contradiction.

$\qedd$\\


\section{Proof of Theorem B}
 We introduce a smooth function
$f:[1,\infty)\to(0,\infty)$ and the half period function for $f.$
Suppose that $f'(1)=0,$ and $f'(x)<0$ on $(1,r_1)$ for some $r_1,$  where $r_1$ is infinite or a finite number greater than 1.
The {\it half period function}
 $\psi_{f}(\nu)$ is defined  for $f$
by
\begin{equation}\label{eq:CT.1}
\psi_{f}(\nu):=2\int_1^{(f|_{(1,r_1)})^{-1}(\nu)}\frac{\nu}{f(x)\sqrt{f(x)^2-\nu^2}}dx
\end{equation}
 for each $\nu\in(f(r_1),f(1)),$ where $f(r_1):=\lim_{x\to\infty}f(x)$ when $r_1$ is infinite.
By imitating the computation in the proof of \cite[Lemma 2.1]{ASTY},
we obtain
\begin{equation}\label{eq:CT.2}
\psi_f(\nu)=-2\int_0^\infty A_f\circ f^{-1}(\sqrt{u(\tau,\nu)})/(a^2\tau^2+1) d\tau,
\end{equation}
where $a:=f(1),$ 
\begin{equation}\label{eq:CT.3}
A_f(x):=\sqrt{a^2-f(x)^2}/f'(x)
\end{equation}
and 
$u(\tau,\nu):=\nu^2(a^2\tau^2+1)/(\tau^2\nu^2+1).$

\begin{lemma}\label{lemCT.1}
Let $f:[1,\infty)\to (0,\infty)$ be the function  defined above.
If 
\begin{equation}\label{eq:CT.4}
f(x)\cdot f'(x)^2+f''(x)(a^2-f(x)^2)\geq0
\end{equation}
holds 
on $(1,r_1),$
then the half period function $\psi_f$ is decreasing on $(f(r_1),f(1)).$
\end{lemma}
\begin{proof}
Since $A_f'(x)=-f'(x)^{-2} (a^2-f(x)^2)^{-1/2}\left(  f(x)\cdot f'(x)^2+f''(x)\cdot (a^2-f(x)^2)\ \right)\leq 0$ by  \eqref{eq:CT.4}, 
the function $A_f$ is decreasing on $(1,r_1).$ 
Hence, for each $\tau>0,$  
 the function $-A_f\circ f^{-1}(\sqrt{u(\tau,\nu)})/(a^2\tau^2+1)$  is decreasing in the variable $\nu.$ 
Therefore,
  the claim  is clear from  \eqref{eq:CT.2}.
  
$\qedd$
\end{proof}
  \begin{remark}
  Note that 
   \eqref{eq:CT.4} holds on $(1,r_1)$ if and only if 
  $A_f(x)^2\cdot(-f''(x)/f(x))\leq 1$ holds on $(1,r_1).$
  \end{remark}

\begin{lemma}\label{lemCT.2}
 If the function ${-f''(x)}/{f(x)}$ is decreasing on $(1,r_1),$
then the half period function $\psi_f$ is decreasing on $(f(r_1),f(1)).$
\end{lemma}
\begin{proof}

From Lemma \ref{lemCT.1} it is sufficient to prove that \eqref{eq:CT.4} holds on $(1,r_1).$
By direct computation, we get
\begin{equation}\label{eq:CT.N1}
F'(x)=(a^2-f(x)^2)\left(  f''(x)/f(x)\right )',
\end{equation}
where $F(x):=f'(x)^2+{f''(x)/f(x)}\cdot (a^2-f(x)^2).$ 
Hence, $F'(x)\geq 0$ on $(1,r_1).$ 
This implies that $F(x)$ is increasing, and $F(x)\geq F(1)=0$ for all $x\in(1,r_1).$
Therefore, by Lemma \ref{lemCT.1}, the half period function $\psi_f$ is decreasing on $(f(r_1),f(1)).$

$\qedd$
\end{proof}

\begin{lemma}\label{lemCT.4}
If there exists a number $r_*\in(1,r_1)$ such that
the function ${-f''(x)}/{f(x)}$ is decreasing on $(1,r_*),$ and $f''(x)\geq 0$
on $[r_*,r_1),$ 
then the half period function $\psi_f$ is decreasing on $(f(r_1),f(1)).$
\end{lemma}
\begin{proof}
From
 \eqref{eq:CT.N1},
$F'(x)\geq 0$ on $(1,r_*).$  This implies that $F(x)$ is increasing on $(1,r_*),$ and 
$F(x)\geq F(1)=0$ for all $x\in(1,r_*).$
Hence \eqref{eq:CT.4} holds on $(1,r_*).$
Since $f''(x)\geq 0, $ and $a^2-f(x)^2\geq 0,$
\eqref {eq:CT.4} also holds for all $x\in(r_*,r_1).$  Therefore, by Lemma \ref{lemCT.1}, the half period function $\psi_f$ is decreasing on $(f(r_1),f(1)).$

$\qedd$
\end{proof}

From now on, we assume that the function $f$ is extended to a smooth function on $[0,\infty)$ such that 
$f'(x)>0 $ on $[0,1),$ $f(0)=0,$ and $f''(1)<0.$
Let $\xi :[0,f(1)]\to [0,1]$ denote the inverse function of $f|_{[0,1]}.$
Since $f'(x)>0$ on $[0,1),$  $\xi$ is smooth on $[0,f(1)),$ and
$\xi'(\nu)=1/f'(\xi(\nu))>0$ and $\xi''(\nu)=-f''(\xi(\nu))/f'(\xi(\nu))^3$ hold 
on $[0,f(1)).$

\begin{lemma}\label{lemCT.5}
$\lim_{x\downarrow 1}(\xi\circ f)'(x)=-1.$
\end{lemma}
\begin{proof}
Since
$$\frac{d}{dx} f'(\xi\circ f(x))^2=2f'(\xi\circ f(x))\cdot f''(\xi\circ f(x))\cdot(\xi\circ f)'(x)
=2f'(x)\cdot f''(\xi\circ f(x)),$$
we get, by l'H\^opital's rule, 
$$\lim_{x\downarrow 1}\frac{f'(x)^2}{f'(\xi\circ f(x))^2}=\lim_{x\downarrow 1}\frac{f''(x)}{f''(\xi\circ f(x))}=1.$$
 Thus, $\lim_{x\downarrow 1}(\xi\circ f)'(x)^2= 
 \lim_{x\downarrow 1}{f'(x)^2}/{f'(\xi\circ f(x))^2}=1.$ This implies that $\lim_{x\downarrow 1}(\xi\circ f)'(x)=-1,$ since 
$(\xi\circ f)'(x)<0$ for all $x\in(1,r_1).$

$\qedd$
\end{proof}

\begin{lemma}\label{lemCT.6}
If the function ${-f''(x)}/{f(x)}$ is decreasing on $(0,r_*)$ for some $r_*\in(1,r_1),$
then 
\begin{equation}\label{eq:CT6n}
1+(\xi\circ f)'(x)\geq 0
\end{equation}
holds 
on $(1,r_*].$
Moreover, \eqref{eq:CT6n} holds on $(1,r_1),$
by adding the assumption that 
 $f''(x)\geq 0$ on $[r_*,r_1).$ 
\end{lemma}

\begin{proof}
Since $-f''(x)/f(x)$ is decreasing on $(0,1]\subset(0,r_*),$ and  since $\xi\circ f(x)\leq 1 $ for all $x\in[1,r_1),$ $-f''(\xi\circ f(x))/f( \xi\circ f(x) )\geq -f''(1)/f(1)>0$ holds for all $x\in[1,r_1).$
In particular, 
\begin{equation}\label{eq:CT.7}
f''(\xi\circ f(x))<0 \: {\rm \:for\: \:all\:\:} x\in[1,r_1).
\end{equation}
It is easy to check that 
\begin{equation}\label{eq:CT.8}
(\xi\circ f)'(x)=\xi'(f(x))f'(x)=\frac{f'(x)}{f'(\xi\circ f(x) )}
\end{equation}
holds on $(1,r_1).$
By differentiating \eqref{eq:CT.8},
it is clear to see that the function $Y(x):=(\xi\circ f)'(x)$ satisfies the following differential equation:

\begin{equation}\label{eq:CT.9}
Y'(x)=\frac{f''(x)-f''(\xi\circ f(x))\cdot Y(x)^2}{f'(\xi\circ f(x))}
\end{equation} 
on $(1,r_1).$
Since the function $-f''(x)/f(x)$ is decreasing on $(0,r_*)$ and since $f(\xi\circ f(x))=f(x),$
$f''(x)\geq f''(\xi\circ f(x))$ holds on $(1,r_*].$
Thus, by \eqref{eq:CT.9},

\begin{equation}\label{eq:CT.10}
Y'(x)\geq\frac{f''(\xi\circ f(x))}{f'(\xi\circ f(x))}\cdot(1- Y(x)^2)
\end{equation}
holds on $(1,r_*].$
Put $M_0:=\sup\{Y(x)^2\: ; \:x\in(1,r_*]  \}.$ By Lemma \ref{lemCT.5}, $M_0\geq1.$
By supposing $M_0>1,$ we will get a contradiction.
Then, $Y(x)^2\leq 1$ implies our first conclusion.
Since $Y(x)^2$ is continuous, 
there exists a number $x_0\in (1,r_*]$ satisfying 
$M_0=Y(x_0)^2$ and $\left.\frac{d}{dx}\right|_{x_0}Y(x)^2\geq0.$ 
By \eqref{eq:CT.10},
$$ 0\leq Y(x_0)Y'(x_0)\leq Y(x_0)\cdot\frac{f''(\xi\circ f(x_0))}{f'(\xi\circ f(x_0))}(1-M_0).$$
Hence, $f''(\xi\circ f(x_0))\geq 0$ which  contradicts \eqref{eq:CT.7}.

Moreover, suppose that $f''(x)\geq 0$  on $[r_*,r_1).$
By \eqref{eq:CT.7} and \eqref{eq:CT.9},  $ (\xi\circ f)''(x)=Y'(x)\geq 0$ on $[r_*,r_1).$ Thus,
$(\xi\circ f)'(x)$ is increasing on $[r_*,r_1).$
Therefore, $1+(\xi\circ f)'(x)\geq1+(\xi\circ f)'(r_*)\geq 0$ for all $x\in(r_*,r_1).$
$\qedd$
\end{proof}

\noindent
{Proof of Theorem B.}\vspace{3mm}
\\
\noindent
We may assume that the number $\lambda=1,$
since the function 
$m_\lambda(x):=\lambda^{-1}m(\lambda x)$
satisfies (M.5) and (M.6),  and $m_\lambda'(1)=0>m_\lambda''(1), \: m_\lambda'(x)>0$ on $[0,1),$ and since
the cut locus structure of the surface of revolution $(\R^2,dr^2+m_\lambda(r)^2d\theta^2)$ is preserved by the homothetic transformation
$f_\lambda: (\R^2,dr^2+m_\lambda(r)^2d\theta^2) \to (\R^2,dr^2+m(r)^2d\theta^2),$ which is defined by
$f_\lambda\circ\exp_o(v):=\exp_o(\lambda v)$ for any tangent vector $v$ at $o,$
where $\exp_o: T_oM\to M$ denotes the exponential map giving a diffeomorphism.


If $m'(x)<0$ on $(1,\infty),$ then, from (M.6) and Lemmas \ref{lemCT.4} and \ref{lemCT.6}, it follows that
$m$ satisfies the properties (M.1), (M.2),  (M.3),  and (M.4). Therefore, by 
{ Theorem A}, the surface $M$ is a  generalized von Mangoldt surface of revolution.

Suppose that  $m'$ has a second zero $r_1(>1).$ Hence, by (M.5), $m'(x)<0$ on $(1,r_1),$ and $m'(r_1)=0.$
If the number $r_1> r_{dc},$
then it is clear  from (M.6)  that $m''(x)\geq 0 $ on $[r_1,\infty),$ and hence the function $m$ satisfies the properties (M.1) and (M.2). By Lemmas \ref{lemCT.4} and \ref{lemCT.6},
$m$ satisfies the properties (M.3) and (M.4).
Suppose that $r_{dc}\geq r_1.$
Since $m'(1)=m'(r_1)=0,$ by applying the mean value theorem for $m',$ there exists a number $r_2\in(1,r_1)\subset(1,r_{dc})$ satisfying $m''(r_2)=0.$
By the property (M.6), $0=(-m''/m)(r_2)\geq (-m''/m)(x)$ for all $x\in[r_2,r_{dc}].$
In particular, $m''(x)\geq0$ on $[r_1,r_{dc}],$ since $[r_1,r_{dc}]\subset [r_2,r_{dc}].$
Therefore, $m$ satisfies (M.1). From Lemmas \ref{lemCT.2}, \ref{lemCT.6}  and (M.6), it is clear that 
$m$ satisfies (M.2), (M.3) and (M.4). Therefore, by {Theorem A}, the surface $M$ is a  generalized von Mangoldt surface of revolution.

$\qedd$


\section{Examples}
We will demonstrate that the function 
$m_0(x):=x/(1+x^2)$ satisfies the properties (M.5) and (M.6).
This example is fundamental, that is, all necessary examples will be constructed from this function.
Since $m_0$ is a smooth  odd function, $m_0'(0)=1$ and $m_0>0$ on $(0,\infty),$  the Riemannian metric $dr^2+m_0(r)^2d\theta^2$  on the Euclidean plane $\R^2,$
where $(r,\theta)$ denotes polar coordinates around the origin $o$ of the plane,
gives  a smooth Riemannian metric of a surface of revolution. It is easy to confirm that
\begin{equation*}\label{eq:3.1}
m_0'(x)=(1-x^2)/(1+x^2)^2,
\quad m_0''(x)=-2x(3-x^2)/(1+x^2)^3
\end{equation*}
and 
\begin{equation}\label{eq:3.2}
(-m_0''/m_0)(x)=2(3-x^2)/(1+x^2)^2.
\end{equation}
Thus, the function satisfies the property   (M.5). 
In our case, the function $m'_0$ does not have a second zero.
By differentiating \eqref{eq:3.2}, we get
\begin{equation}\label{eq:3.3}
(-m_0''/m_0)'(x)={4x(x^2-7)}/{(1+x^2)^3}.
\end{equation}
By \eqref{eq:3.3},
it is clear that 
the  function $(-m_0''/m_0)(x)$ is decreasing on $(0,\sqrt{7}).$
Therefore, the function $m_0$ satisfies the properties  (M.5) and (M.6).
By { Theorem B},
$(\R^2,dr^2+m_0(r)^2d\theta^2)$ is a (non-trivial) generalized von Mangoldt surface of revolution such that $\lim_{x\to\infty}m'_0(x)=0,$ and $\int_0^\infty|m_0''(x)|dx<\infty.$ 
Summing up the argument above, we get:
\begin{proposition}\label{prop3.1}
The surface 
$(\R^2,dr^2+m_0(r)^2d\theta^2)$ 
is  a generalized von Mangoldt surface of revolution such that 
the function
 $m_0$
satisfies  the properties (M.5) and  (M.6) with $\lambda=1,r_{dc}=\sqrt{7},$ and furthermore,
$\lim_{x\to\infty}m'_0(x)=0,$ $\int_0^\infty|m''_0(x)|dx<\infty$ and $m_0''(x)>0$ on $[\sqrt{7},\infty).$

\end{proposition}


The main aim of this section is to prove the following proposition.
\begin{proposition}\label{prop3.3}
For each  number $\alpha>0,$
there exists a generalized von Mangoldt surface of revolution
$(\R^2,dr^2+m_\alpha(r)^2d\theta^2)$ such that the function
 $m_\alpha$
satisfies  the properties (M.5) and  (M.6) with $\lambda=1,r_{dc}=\sqrt{7},$ and furthermore,
$\lim_{x\to\infty}m'_\alpha(x)=\alpha,$ $\int_0^\infty|m''_\alpha(x)|dx<\infty$ and $m_\alpha''(x)>0$ on $[\sqrt{7},\infty).$ 
\end{proposition}

\begin{proof}

Choose   any positive number $\alpha.$
We define a smooth function  
$m_{\alpha}$ by
\begin{equation}\label{eq:3.9}
m_{\alpha} (x):=m_0(x)+\alpha\int_0^x\phi_1(x)dx,
\end{equation}
where $m_0(x)=x/(1+x^2),$  and $\phi_1(x)$ denotes a smooth step function, i.e.,  
$\phi_1(x)=0$ for $x\leq r_{dc}:=\sqrt{7}, \phi_1(x)=1$ for $x\geq r_{dc}+1,$ and $\phi_1$ is increasing on $\R.$
Since $m_\alpha(x)=m_0(x)$ on $[0,r_{dc}],$ and since $m_\alpha''(x)=m_0''(x)+\alpha\phi_1'(x)\geq m_0''(x)$ on $[0,\infty),$ it is clear that $m_\alpha$ satisfies the properties (M.5) and (M.6),
and furthermore, it is clear from \eqref{eq:3.9} that 
$\lim_{x\to\infty}m_\alpha'(x)=\alpha,$  $\int_0^\infty|m_\alpha''(x)|dx<\infty$ and  $m_\alpha''(x)>0$ on $[\sqrt{7},\infty).$

$\qedd$
\end{proof}

\section{Oscillations of the Gaussian curvature}
In this section, {Theorem C} in the introduction will be proven. Some lemmas and one proposition are necessary for this.
\begin{definition}
For each positive integer $n$, we set 
$I_n := [-1 + 2n,1 + 2n]$ and 
define two norms, $|\cdot|_{1,n}$
and $|\cdot|_{2,n}$,
for each smooth function $f:I_n\to\R$ by
\begin{gather*}
  |f|_{1,n} := \sup\{|f(x)| + |f'(x)|;x 
  \in I_n\} \quad \mbox{and}\\
  |f|_{2,n} := \sup\{|f(x)| + |f'(x)| + |f''(x)|;x \in I_n\}, 
\end{gather*} 
respectively.
\end{definition}

\begin{lemma}\label{lem4.1}
For each $\epsilon>0,$ there exists a smooth function $\tilde f_\epsilon : \R\to \R$
satisfying the following \eqref{eq:4.1}, \eqref{eq:4.2}, and \eqref{eq:4.3}:
\begin{equation}\label{eq:4.1}
\supp \tilde f_\epsilon\subset (-1,1).
\end{equation}
\begin{equation}\label{eq:4.2}
\max\{|\tilde f_\epsilon(x)|+|\tilde f_\epsilon'(x)|+|\tilde f_\epsilon''(x)|\: ; \: x\in \R\}<\epsilon.
\end{equation}
\begin{equation}\label{eq:4.3}
\max \{\tilde f_\epsilon'''(x); x\in \R \}>\frac{1}{\epsilon} \quad and \quad \min \{\tilde f_\epsilon'''(x); x\in \R\} <\frac{-1}{\epsilon}.
\end{equation}

\end{lemma}
\begin{proof}
Choose any $\epsilon>0$ and  a positive integer $k$ satisfying 
$$K:=\sqrt{(4k+1)\pi}>\max \{3+3|\phi'(x)|+|\phi''(x) |  \: ; \: x\in \R \}/\epsilon,$$
where $\phi : \R\to \R$ denotes a smooth function with $\supp \phi\subset(-1,1),$
satisfying $\phi=1$ on $[-1/2,1/2].$
Define a smooth function $\tilde f_\epsilon(x)$
by $\tilde f_\epsilon(x):=\phi(x)\cos(K^2x)/K^5.$

By the triangle inequality,
we obtain
\begin{equation*}
|\tilde f_\epsilon'(x)|\leq |\phi'(x)/K^5|+1/K^3\leq(1+|\phi'(x)|)/K,
\quad 
|\tilde f_\epsilon''(x)|\leq( |\phi''(x)|+2|\phi'(x)|+1)/K.
\end{equation*}
Hence,
\begin{equation*}
|\tilde f_\epsilon(x)|+|\tilde f_\epsilon'(x)|+|\tilde f_\epsilon''(x)|\leq
(3+3|\phi'(x)|+|\phi''(x)|)/K<\epsilon.
\end{equation*}
Then, it is clear that  the function $\tilde f_\epsilon(x)$
satisfies \eqref{eq:4.1} and \eqref{eq:4.2}. Let us check the property \eqref{eq:4.3}.
Since $\phi(\pm1/2)=1, \phi'(\pm1/2)=\phi''(\pm1/2)=\phi'''(\pm1/2)=0,$
we obtain
$$\tilde f_\epsilon'''(\pm1/2)=K\sin(K^2(\pm1/2))=\pm K.$$
Hence, \eqref{eq:4.3} clearly  holds for $\tilde f_\epsilon.$

$\qedd$
\end{proof}

\begin{lemma}\label{lem4.2}
Let $\{\epsilon_n\}_{n\geq n_0}$ be any sequence of positive numbers, where 
$n_0$ denotes  a positive integer.
Then for each $n\geq n_0$ there exists a smooth function $f_n : \R\to \R$
such that

\begin{equation}\label{eq:4.4}
\supp f_n\subset(-1+2n,1+2n),
\end{equation}
\begin{equation}\label{eq:4.5}
|f_n|_{2,n}<\epsilon_n,
\end{equation}
\begin{equation}\label{eq:4.6}
\max \{ f_{n}'''(x); x\in I_n\}>\frac{1}{\epsilon_n} \quad and \quad \min\{ f_{n}'''(x); x\in I_n\}<\frac{-1}{\epsilon_n}.
\end{equation}

\end{lemma}
\begin{proof}
By applying Lemma \ref{lem4.1}, for each $\epsilon_n$ there exists a smooth function $\tilde f_{\epsilon_n}$ satisfying \eqref{eq:4.1}, \eqref{eq:4.2} and \eqref{eq:4.3}
with $\epsilon=\epsilon_n.$
Then the function $f_n(x):=\tilde f_{\epsilon_n}(x-2n)$ satisfies 
\eqref{eq:4.4},    
\eqref{eq:4.5} and \eqref{eq:4.6} for each $n.$
$\qedd$

\end{proof}

\begin{definition}
A sequence $\{f_n\}$ of smooth functions  is said to be {\it associated with} 
a   sequence $\{\epsilon_n\} $ of positive numbers 
if the functions $f_n:\R\to\R$ satisfy \eqref{eq:4.4}, \eqref{eq:4.5} and \eqref{eq:4.6}
in Lemma \ref{lem4.2}.
\end{definition}

\begin{proposition}\label{prop4.4}
For any smooth function $\hat m: (0,\infty)\to (0,\infty)$ 
satisfying 
\begin{equation*}\label{eq:4.7}
\hat m''(x)>0 {\quad on} \quad [t_0,\infty)\quad \mbox{for some} \quad t_0>0,
\end{equation*}
and
\begin{equation*}\label{eq:4.8}
\hat m'(\infty):=\lim_{x\to\infty}\hat m'(x)\geq 0,
\end{equation*}
there exists a sequence $\{C_n\}_{n\geq n_0}$ of positive numbers (depending only on the function $\hat m$) such that
for any sequence $\{f_n\}_{n\geq n_0}$ associated with the sequence $\{C_n\}_{n\geq n_0},$
the function $m(x):=\hat m(x)(1+\sum_{n=n_0}^\infty f_n(x))$
satisfies

\begin{equation}\label{eq:4.9}
m(x)=\hat m(x)\quad on\quad (0,t_0],
\end{equation}
\vspace*{-0.8ex}
\vspace*{-0.8ex}
\begin{equation}\label{eq:4.10}
m''(x)>0\quad on \quad [t_0,\infty),
\end{equation}
\begin{equation}\label{eq:4.11}
\lim_{x\to\infty}m'(x)=\hat m'(\infty),
\end{equation}
\vspace*{-0.8ex}
\begin{equation}\label{eq:4.12}
m''(x)/m(x) \quad \mbox{is not monotone on} \quad [a,\infty)\quad \mbox{for any} \quad a>0.
\end{equation}
Moreover, if $\int_0^\infty|\hat m''(x)|dx<\infty,$ then $\int_0^\infty|m''(x)|dx<\infty.$

\end{proposition}

\begin{proof}
Choose any integer  $n_0$ greater than $(1+t_0)/2,$ and fix it.
For each integer $n\geq n_0,$
we define a positive number $C_n$ by 
\begin{equation} \label{eq:4.13}
  C_n := \min \big\{a_n(2|\hat m|_{2,n})^{-1}, 
  (3[1 + 4|\hat m'/\hat m|_{1,n} + |\hat m''/\hat m|_{1,n}])^{-1}, 
      (n^2|\hat m|_{2,n})^{-1}\big\},
\end{equation} 
where $a_n := \min \{\hat m''(x) ; x \in I_n \}.$

Let  $\{f_n\}_{n\geq n_0}$ be any sequence of smooth functions associated with the sequence 
$\{ C_n\}_{n\geq n_0}$.
Let us check  \eqref{eq:4.9} first.
Since $\supp f_n\subset I_n$ 
and $I_n$ is disjoint from $(0,t_0]$ for each $n\geq n_0,$
it is clear that $m(x)=\hat m(x)$ holds on $(0,t_0].$

By direct computation,
we get
$$ \hat m''(x)-m''(x)=-(f_n''\hat m+2f_n'\hat m'+f_n\hat m'')(x)\:\:\mbox{on} \:\: I_n.$$
By the triangle inequality and \eqref{eq:4.13},
\begin{equation}\label{eq:4.14}
|\hat m''(x)-m''(x)|\leq 2|f_n|_{2,n}\cdot|\hat m|_{2,n}< 2C_n|\hat m|_{2,n}\leq a_n\:\: \:\mbox{on}\:\: I_n.
\end{equation} 
Thus,  \eqref{eq:4.10} is clear from
 \eqref{eq:4.14}.

Since $m'(2k+1)=\hat m'(2k+1)$ holds for any positive integer $k>t_0 ,$
$\lim_{k\to\infty}m'(2k+1)=\lim_{k\to\infty}\hat m'(2k+1)=\hat m'(\infty)$ holds.
On the other hand,
 $m'$ is increasing on $[t_0,\infty)$ by \eqref{eq:4.10}. Hence,  \eqref{eq:4.11} holds.

Since $|f_n(x)|+|f_n'(x)|+|f_n''(x)|\leq|f_n|_{2,n}<C_n\leq 1/3<1/2$ holds on $I_n,$
the inequalities 
$1/2<1+f_n(x)<3/2,$ \:
$|(1/(1+f_n))'(x)|<2,$\:  $|f'_n/(1+f_n)(x)|<1$ and 
$|(f_n'/(1+f_n))'(x)|<2$ 
hold on $I_n.$

By direct computation, we obtain
\begin{equation}\label{eq:4.15}
(m''/m)(x)=f_n''/(1+f_n)(x)+2f_n'/(1+f_n)(x)\cdot(\hat m'/\hat m)(x)+(\hat m''/\hat m)(x)\: \:\:\mbox{on}\:\: I_n.
\end{equation}
By differentiating \eqref{eq:4.15} and  by the triangle inequality, we get
\begin{equation*}\label{eq:4.16}
|(m''/m)'(x)-f_n'''/(1+f_n)(x)|\leq1+4|\hat m'/\hat m|_{1,n}+|\hat m''/\hat m|_{1,n}\:\:\mbox{on}\: I_n.
\end{equation*}
Thus, by \eqref{eq:4.13},
we obtain
\begin{equation*}
{-1}/{(3C_n)}\leq(m''/m)'(x)-f_n'''/(1+f_n)(x)\leq 1/(3C_n){\quad\rm on\quad} I_n,
\end{equation*}
so that it is clear that
\begin{equation*}
\max_{x\in I_n}(m''/m)'(x)- f_n'''/(1+f_n)(x)\geq -1/(3C_n)
{\quad\rm on\quad} I_n,
\end{equation*}
and
\begin{equation*} 
\min_{x\in I_n}(m''/m)'(x)- f_n'''/(1+f_n)(x)\leq 1/(3C_n)
{\quad\rm on\quad} I_n.
\end{equation*}
In particular,
\begin{equation}\label{eq:4.17}
\max_{x\in I_n}(m''/m)'(x)- \max_{x\in I_n}f_n'''/(1+f_n)(x)\geq -1/(3C_n)
\end{equation}
and
\begin{equation}\label{eq:4.18}
\min_{x\in I_n}(m''/m)'(x)- \min_{x\in I_n}f_n'''/(1+f_n)(x)\leq 1/(3C_n).
\end{equation}
Since  $\max_{x\in I_n}f_n'''(x)>1/C_n,$  $\min_{x\in I_n}
 f_n'''(x)<-1/C_n, $  and  since $1/2<1+f_n(x)<3/2$ holds on $I_n,$
 we get
 \begin{equation*}
 \max_{x\in I_n}f_n'''/(1+f_n)(x)>{2}/{(3C_n)},\quad \min_{x\in I_n}f_n'''/(1+f_n)(x)< -2/(3C_n),
 \end{equation*}
 and hence by \eqref{eq:4.17} and \eqref{eq:4.18},
$\max_{x\in I_n}(m''/m)'(x)>{2}/{(3C_n)}-{1}/{(3C_n)}>0$ and
$\min_{x\in I_n}(m''/m)'(x)<{-2}/{(3C_n)}+{1}/{(3C_n)}<0.$ 
This implies that \eqref{eq:4.12} holds for 
the function $m.$ 
Suppose that  $\int_0^\infty|\hat m''(x)|dx<\infty.$
By \eqref{eq:4.13} and \eqref{eq:4.14},
$|m''(x)|\leq|\hat m''(x)|+|m''(x)-\hat m''(x)|<|\hat m''(x)|+2C_n|\hat m|_{2,n}\leq |\hat m''(x)|+2/n^{2}$ holds 
on $I_n.$ Thus, $\int_0^\infty |m''(x)|dx<\infty$ if $\int_0^\infty|\hat m''(x)|dx<\infty.$

$\qedd$
\end{proof}

\noindent
{Proof of Theorem C.  }\vspace{3mm}\\
Let $\widetilde M=(\R^2,dr^2+\tilde m(r)^2d\theta^2)$
denote a surface of revolution with finite total curvature $c.$
Since $\widetilde M$ admits a finite total curvature, we get
$\int_0^\infty|\tilde m''(x)| dx<\infty.$ 
Thus
\begin{equation}\label{eq:5.18}
c=-\int^\infty_0  
\tilde m''(r) dr \int^{2\pi}_0 d\theta
=2\pi(1-\lim_{ r\to\infty}\tilde m'(r)).
\end{equation}
If we  apply {Theorem B}, {Propositions \ref{prop3.1}, \ref{prop3.3},  and \ref{prop4.4},
we obtain a generalized von Mangoldt surface of revolution $M:=(\R^2,dr^2+m(r)^2d\theta^2)$ with $\lim_{r\to\infty}m'(r)=(2\pi-c)/(2\pi)$ and with 
$\int_0^\infty |m''(x)|dx<\infty$
such that
the Gaussian curvature along each meridian is not monotone on $[a,\infty)$
for any $a>0.$
By \eqref{eq:5.18},
it is clear that the total curvature of $M$ is $c.$
$\qedd$

\begin{remark}
Note that if the surface $\widetilde M$ admits a total curvature in the 
sense of Cohn-Vossen (see \cite[Definition 2.1.3]{SST}), then
 $c\leq 2\pi$ by 
a theorem proved by Cohn-Vossen (see \cite[Theorem 2.1.1]{SST}).
 Note also that the total curvature of $\widetilde M$ is finite if and only if $\int_0^\infty|\tilde m''(x)|dx<\infty.$
\end{remark}
\begin{remark}
{Theorem C} is still true for any surface of revolution admitting infinite total curvature.
However, such surfaces of revolution are not useful in radial curvature geometry.
In fact, surfaces of revolution used as model surfaces are always assumed to admit a finite total curvature. For examples, see \cite{KT1,KT2,KT3,KT4,TK}.
\end{remark}

\section{Asymmetric 2-spheres of revolution with simple cut locus structure}

By imitating  some techniques used in the proof of { Theorem} A, we can prove the existence of an  {\it asymmetric} 2-sphere of revolution with simple cut locus structure.

In order to define  an {\it asymmetric 2-sphere of revolution,}
let us first  review  a 2-sphere of revolution $(\Sph^2,g).$ We start from the 2-dimensional unit sphere $(\Sph^2,dr^2+\sin^2 rd\theta^2),$  where $(r,\theta)$ denotes geodesic polar coordinates around the north pole $N.$
In order to introduce a new metric $g$ to the unit sphere, we need a smooth function
$m: [0,\pi]\to [0,\infty)$ satisfying 
\begin{equation*}\label{eq:As.1}
m(0)=m(\pi)=0,\quad m'(0)=-m'(\pi)=1,
\end{equation*}
and
\begin{equation*}\label{eq:As.2}
m(r)>0 \quad {\rm for \: all\:} r\in(0,\pi).
\end{equation*}
Then one can introduce a smooth Riemannian metric 
\begin{equation}\label{eq:As.3}
g=dr^2+m(r)^2d\theta^2
\end{equation}
on the open subset $\Sph^2\setminus\{N,S\}$ of the unit sphere, where $S$ denotes the unique cut point of the point $N,$ which is called the {\it south pole.}
If the functions $m(r)$ and $m(\pi-r)$ are extendable to a smooth odd function around $r=0$ respectively, then it was proved in \cite[Theorem 7.1.1]{SST} that the Riemannian 
metric $g$ is extendable to a smooth one on the entirety of $\Sph^2.$

In this article, the Riemannian manifold $(\Sph^2,g)$ with smooth Riemannian metric $g$ in \eqref{eq:As.3} is called a {\it 2-sphere of revolution}.
In general, the function $m$ does not satisfy
$m(r)=m(\pi-r),$ i.e., the 2-sphere of revolution $(\Sph^2,dr^2+m(r)^2d\theta^2)$ does not need to be reflectively symmetric with respect to the parallel $r=\pi/2.$
A 2-sphere of revolution is said to be {\it asymmetric} if the sphere does not have a reflective isometry with respect to the parallel $r=\pi/2.$

In this section, {\it asymmetric 2-spheres of revolution with simple cut locus structure
will be introduced.}
Before describing such an example, we need the following theorem which gives a sufficient condition for a 2-sphere   of revolution to have a simple cut locus structure.
\begin{theorem}\label{thAs.1}
If a 2-sphere of revolution $(\Sph^2,dr^2+m(r)^2d\theta^2)$ satisfies the following three properties, then, for each point $p\in\theta^{-1}(0),$ the cut locus of the point $p$ is a single point or a subarc of $\theta^{-1}(\pi).$
\begin{description}
\item (A.1)\quad
The derivative  $m'(x)$ of the function $m(x)$ has a unique zero $r_e\in(0,\pi).$  
\item
 (A.2)\quad  The function $-m''(x)/m(x)$ is  decreasing  on $ (0,r_e) $
and increasing on $(r_e,\pi).$
\item
(A.3)\quad  The function $m(x)$ satisfies either 
\begin{equation*}
\quad\quad 1+(\xi\circ m)'(x)\geq 0 \quad on \quad (r_e,\pi),\quad or \quad
 1+(\xi\circ m)'(x)\leq 0 \quad on \quad (r_e,\pi),
\end{equation*}
\quad \: where $\xi(\nu):=(m|_{[0,r_e]})^{-1}(\nu)$ on $[0,m(r_e)].$

\end{description}
\end{theorem}
\begin{remark}\label{remA.2}
If the 2-sphere of revolution in Theorem \ref{thAs.1} admits a reflective isometry with respect to the parallel $r=\pi/2,$ then $r_e=\pi/2$ and the property (A.3) is satisfied, since $\xi\circ m(r)=\pi-r$ holds on $[\pi/2,\pi].$ 
\end{remark}
\begin{remark}\label{remA.3}
Since the Gaussian curvature of a 2-sphere $(\Sph^2,dr^2+m(r)^2d\theta^2)$  of revolution equals $(-m''/m)(r(q))$ at each point $q\in r^{-1}( (0,\pi) ),$
the property (A.2) implies that the Gaussian curvature is decreasing along each meridian 
from  the north pole (respectively the south pole) to the point on the parallel $r=r_e.$
\end{remark}
\begin{proof}
Since 
$r_e$ is the unique zero of $m'$ by the property (A.1), the parallel $r=r_e$ is a geodesic and no other parallel is a geodesic. In our case, the corresponding half period functions $\varphi(\nu),\psi(\nu)$  on $(0,m(r_e))$ are defined by
$$\varphi(\nu):=2\int^{r_e}_{\xi(\nu)} f(x,\nu)dx\quad\mbox{and}\quad\psi(\nu):=2\int_{r_e}^{\eta(\nu)} f(x,\nu)dx,$$
where $f(x,\nu):=\nu/(m(x)\sqrt{m(x)^2-\nu^2}), \: \xi(\nu):=\min m^{-1}(\nu)$ and $\eta(\nu):=\max m^{-1}(\nu).$
From the property (A.2), $\varphi(\nu)$ and $\psi(\nu)$ are increasing on $(0,m(r_e)),$ respectively, and, in particular, for all $\nu>0,$
$\varphi(\nu)\geq\lim_{\nu\downarrow 0}\varphi(\nu)=\pi$ and $\psi(\nu)\geq\lim_{\nu\downarrow 0}\psi(\nu)=\pi.$
Therefore, the subarc $r^{-1}(r_e)\cap \theta^{-1}([0,\pi])$ of the parallel $r=r_e$ is
a minimal geodesic segment.

Suppose that there exists  a point $p\in\theta^{-1}(0)$ having a cut point $q\notin\theta^{-1}(\pi).$
We may assume that $q\in \theta^{-1}( (0,\pi)  ) $ and that $q$ is a conjugate point of $p$
along a minimal geodesic segment $\alpha$ joining $p$ to $q.$
Otherwise, there exists another  minimal geodesic segment $\beta$ joining $p$ to $q.$
Then, by imitating the same method used in the proof in \cite[Lemma 3.1]{ST2}, we can find a cut point $q_0$ of $p$ in the domain bounded by $\alpha$ and $\beta$ which is a conjugate point of $p$ along any minimal geodesic segments joining $p$ to $q_0.$
Therefore, there exists a cut point $q\in\theta^{-1}( (0,\pi) )$ of $p$ which is a conjugate point of $p$ along a minimal geodesic segment $\alpha$ joining $p$ to $q.$
Since $\alpha$ is not  a subarc of the parallel $r=r_e,$ and since  $q\in\theta^{-1}( (0,\pi) ),$ the Clairaut constant $\nu_0$ of $\alpha$ is positive and  less than $m(r_e).$ 
From the property (A.2), and by imitating the proof of Lemma 3.1 in \cite{ST2}, we may assume that $r(p)\in(r_e,\pi)$ and $r(q)\in(0,r_e).$ Moreover, we may assume that $\alpha=\gamma_{\nu_0}^{(p)}.$ Recall that $\gamma_{\nu_0}^{(p)}$ denotes the geodesic with Clairaut constant $\nu_0$ emanating from $p=\gamma_{\nu_0}^{(p)}(0),$ which is defined in Definition \ref{def2.3}.
Hence the function $\Phi_{\nu_0}(t,h(t))$ is less than $\pi$ for some $t\in(r_e,\eta(\nu_0)],$
where $\Phi_{\nu_0}(u,v):=\varphi(\nu_0)+\int ^u_{r_e}f(x,\nu_0)dx-\int^{r_e}_v f(x,\nu_0)dx$
on $[r_e,\eta(\nu_0) )\times(\xi(\nu_0),r_e],$
and $h(t)$ is the function defined by the equation
$\varphi'(\nu_0)+\int^t_{r_e}f_{\nu}(x,\nu_0)dx-\int^{r_e}_{h(t)} f_{\nu}(x,\nu_0)dx=0,$ and
$f_{\nu}(x,\nu):=m(x)/\sqrt{m(x)^2-\nu^2}^{3}$
in this case.
Since $\Phi_{\nu_0}(r_e,h(r_e))\geq\pi$ and $\lim_{t\uparrow\eta(\nu_0)}\Phi_{\nu_0}(t,h(t))=(\varphi(\nu_0)+\psi(\nu_0))/2\geq\pi,$
the function $\Phi_{\nu_0}(t,h(t))$ attains a minimum  at some $t=t_0\in(r_e,\eta(\nu_0))$ and 
$\Phi_{\nu_0}(t_0,h(t_0))<\pi.$ 
Thus, $\left.\frac{d}{dt}\right|_{t_0}\Phi_{\nu_0}(t,h(t))=0,$ and we obtain $h(t_0)=\xi\circ m(t_0).$ By this equation $h(t_0)=\xi\circ m(t_0)$ and  the following lemma,
$\pi>\Phi_{\nu_0}(t_0,h(t_0)) =\Phi_{\nu_0}(t_0,\xi\circ m(t_0)) \geq\pi.$ 
This  is a contradiction.
$\qedd$
\end{proof}
\begin{lemma}\label{lemA.2}
For any 
$t\in [r_e,\eta(\nu_0)],$
$\Phi_{\nu_0}(t,\xi\circ m(t))\geq\pi.$
\end{lemma}
\begin{proof}
Since $\frac{d}{dt}\Phi_{\nu_0}(t,\xi\circ m(t))=f(t,\nu_0)(1+(\xi\circ m)'(t)),$
the function $\Phi_{\nu_0}(t,\xi\circ m(t))$ is either increasing or decreasing on $[r_e,\eta(\nu_0)]$ by the property (A.3).
It is clear that $\Phi_{\nu_0}(r_e,\xi\circ m(r_e) )=\varphi(\nu_0)\geq\pi,$ and $\Phi_{\nu_0}(\eta(\nu_0),\xi\circ m(\eta(\nu_0)) )=(\varphi(\nu_0)+\psi(\nu_0))/2\geq\pi.$
Hence, $\Phi_{\nu_0}(t,\xi\circ m(t))\geq \pi$ for all $t\in [r_e,\eta(\nu_0)].$

$\qedd$
\end{proof}

Let $\{M_\lambda\}_\lambda:=\{(\Sph^2,dr^2+m_\lambda(r)^2d\theta^2);\lambda>-1\}$ denote a family of 2-spheres $(\Sph^2,dr^2+m_\lambda(r)^2d\theta^2)$ of revolution obtained by rotating a family of smooth Jordan curves \\ $(l_\lambda\cdot x_\lambda(u),l_\lambda\cdot y(u),0),$
$u\in[0,2\pi],$ around the $x$-axis in Euclidean space.
Here,
\begin{equation}\label{eq:As.4}
x_\lambda(u):=\cos u+\lambda(\cos u-(\cos^4u)/12) 
\quad\mbox{and}\quad y(u):=\alpha\sin u,
\end{equation}
where  $l_\lambda:=\pi(\int^\pi_0\sqrt{x_\lambda'(u)^2+y'(u)^2} du)^{-1},$ and $\alpha\in(0,1)$ is a fixed constant.
Note that for each $\lambda>-1,$ the plane curve $(l_\lambda x_\lambda(u),l_\lambda y(u),0) $ is regular at each $u\in[0,2\pi],$ i.e., its velocity vector is nonzero at each $u\in[0,2\pi],$  its length  equals $2\pi,$
and, for  each $\lambda\in(-1,0)\cup(0,\infty),$ the 2-sphere $M_\lambda$ is asymmetric.

By making use of classical differential geometry, we obtain
\begin{equation*}\label{eq:As.5}
m_\lambda(r)=l_\lambda y(u(r)),
\end{equation*}
where $u(r)$ denotes the inverse function of the function $r=X(u),u\in \R,$
which is defined by
$X(u):=l_\lambda\int_0^u \sqrt{x_\lambda'(u)^2+y'(u)^2} du.$
Hence, 
\begin{equation}\label{eq:As.6}
\frac{du}{dr}=\frac{1}{l_\lambda\sqrt{x_\lambda'(u(r))^2+y'(u(r))^2} }
\end{equation}
and 
\begin{equation}\label{eq:As.7}
m_\lambda'(r)=\frac{y'(u(r))} {\sqrt{x_\lambda'(u(r))^2+y'(u(r))^2} }.
\end{equation}
From  direct computation, it follows that
$$m_\lambda''(r)=\frac{x_\lambda'(u(r))\cdot[ x_\lambda'(u(r))y''(u(r))-x_\lambda''(u(r))y'(u(r)) ]}{l_\lambda(x_\lambda'(u(r))^2+y'(u(r))^2)^2}$$
and
\begin{equation}\label{eq:As.8}
\frac{-m_\lambda''(r)}{m_\lambda(r)}=\frac{Q_\lambda (u(r))}{l_\lambda^2(x_\lambda'(u(r))^2+y'(u(r))^2)^2},
\end{equation}
where
$Q_\lambda(u):=-x'_\lambda(u)\cdot(x_\lambda'(u)y''(u)-x_\lambda''(u)y'(u)) \cdot y(u)^{-1}.$
Note that  \eqref{eq:As.8} is also a well-known result in classical differential geometry.

By \eqref{eq:As.4}, we obtain
\begin{equation}\label{eq:As.9}
x_\lambda'(u)=-\sin u\cdot(1+F_\lambda'(\cos u)),
\end{equation}
where 
$$F_\lambda(x):=\lambda(x-x^4/12).$$
Note here that \begin{equation}\label{eq:As.9n}
|F_\lambda'(\cos u)|\leq 4|\lambda|/3
\end{equation}
holds for all $u\in[0,\pi].$
Since $x_\lambda''(u)=-\cos u\cdot(1+F_\lambda'(\cos u))+\sin ^2u\cdot F_\lambda''(\cos u),$ $y'(u)=\alpha\cos u$ and $ y''(u)=-\alpha\sin u,$ we get
$Q_\lambda(u)=(1+F_\lambda'(\cos u))^2-\cos u\cdot\sin ^2u\cdot (1+F_\lambda'(\cos u))\cdot F_\lambda''(\cos u).$
Since $F_\lambda''(x)=-\lambda x^2$ and $\sin 2u=2\cos u\cdot \sin u,$
we get
\begin{equation}\label{eq:As.10}
Q_\lambda(u)=(1+F_\lambda'(\cos u))^2+\lambda\cos u\cdot\sin^22u\cdot(1+F_\lambda'(\cos u))/4.
\end{equation}

From now on, we will prove that for any small $|\lambda|,$ the function
$m_\lambda$ satisfies the properties (A.1), (A.2), and (A.3). Therefore,
{\it if $|\lambda|>0$ is sufficiently small, by Theorem \ref{thAs.1}, for each point $p\in\theta^{-1}(0)$ on the asymmetric  
2-sphere $(\Sph^2,dr^2+m_\lambda(r)^2d\theta^2),$   
the cut locus of $p$ is a single point or a subarc of $\theta^{-1}(\pi).$}

\begin{lemma}\label{lemAs.5}
The number $r_e:=u^{-1}(\pi/2)\in(0,\pi)$ is a unique zero of the derivative  $m_\lambda'$ of $m_\lambda.$
\end{lemma}
\begin{proof}
From \eqref{eq:As.4} and \eqref{eq:As.7}, it is clear that    $m_\lambda'(r_e)=0.$  Since the function $u(r)$ is strictly increasing on $[0,\pi],$ the number $r_e$  is a unique zero of  $m_\lambda'|_{[0,\pi]}.$
\end{proof}
$\qedd$
\begin{lemma}\label{lemAs.6}
For all $\lambda\geq 0$ {\rm(}respectively 
$\lambda\in(-1,0]
{\rm)},$
$1+(\xi\circ m_\lambda)'(x)\geq 0$ {\rm (}respectively $1+(\xi\circ m_\lambda)'(x)\leq 0${\rm)} holds on $(r_e,\pi).$
\end{lemma}
\begin{proof}
Choose any $\lambda\geq 0$ (respectively $\lambda\in(-1,0]$) and fix it.
It is sufficient to prove that $m_\lambda'(r_1)+m_\lambda'(r_2)\geq 0$ (respectively
$m_\lambda'(r_1)+m_\lambda'(r_2)\leq 0$) holds for any $r_1\in(r_e,\pi),$  where $r_2:=\xi\circ m_\lambda(r_1),$ since
$1+(\xi\circ m_\lambda)'(r_1)=1+\xi'( m_\lambda(r_2)  )\cdot m_\lambda'(r_1)=
1+ m_\lambda '  (r_1) / m_\lambda'(r_2)=(m_\lambda'(r_2)+m_\lambda'(r_1))/m_\lambda'(r_2).$
Note that the equation $\xi\circ m_\lambda(r_1)=r_2$ implies that
\begin{equation}\label{eq:6.11}
l_\lambda\sin u(r_1)=\alpha^{-1}m_\lambda(r_1)=\alpha^{-1}m_\lambda(r_2)=l_\lambda\sin u(r_2).
\end{equation}
Since the function $u(r)$ is strictly increasing and $0<r_2<r_e<r_1<\pi,$ we get, by \eqref{eq:6.11},
\begin{equation}\label{eq:As.12}
0<u_2:=u(r_2)<\pi/2=u(r_e)<u_1:=u(r_1)<\pi,
{\:\rm and\: }
u_1+u_2=\pi.
\end{equation}
By \eqref{eq:As.4}, \eqref{eq:As.12}, and 
 \eqref{eq:As.7}, we get
\begin{equation*}\label{eq:As.13}
y(u_1)=y(u_2) ,\quad \cos u_1=-\cos u_2<0
\end{equation*}
and
\begin{equation}\label{eq:As.14}
m'_\lambda(r_2)+m'_\lambda(r_1)=\alpha\cos u_2\left( \frac{1}{\sqrt{x'_\lambda(u_2)^2+y'(u_2)^2}}-\frac{1}{\sqrt{x'_\lambda(u_1)^2+y'(u_2)^2} } \right).
\end{equation}
Since $F_\lambda'(x)=\lambda(1-x^3/3),$
 we obtain
\begin{equation*}\label{eq:As.15}
F_\lambda'(x)+F_\lambda'(-x)=2\lambda,\quad 
F_\lambda'(x)-F_\lambda'(-x)=-2\lambda x^3/3.
\end{equation*}
By \eqref
{eq:As.9},
we obtain
\begin{equation*}\label{eq:As.16}
x_\lambda'(u_2)=-\sin u_2\cdot(1+F_\lambda'(\cos u_2)),\quad 
x_\lambda'(u_1)=-\sin u_2\cdot(1+F_\lambda'(-\cos u_2)).
\end{equation*}
Hence,
\begin{equation}\label{eq:As.17}
x_\lambda'(u_1)^2-x_\lambda'(u_2)^2=4/3\cdot\lambda(1+\lambda)\cdot\sin^2u_2\cdot\cos ^3u_2.
\end{equation}
Now the claim of our lemma is clear from \eqref{eq:As.14} and \eqref{eq:As.17}.
$\qedd$
\end{proof}

From \eqref{eq:As.8}, it follows that 
\begin{equation}\label{eq:As.18}
\frac{d}{dr}\:\frac{-m''_\lambda(r)  }{m_\lambda(r)  }=\frac{P_\lambda(u(r))}{l_\lambda^2\cdot(x_\lambda'(u(r))^2+y'(u(r))^2)^3}\cdot\frac{du}{dr},
\end{equation}
where 
\begin{equation}\label{eq:As.19}
P_\lambda(u):=Q'_\lambda(u)(x_\lambda'(u)^2+y'(u)^2)-2Q_\lambda(u)
(x_\lambda'(u)^2+y'(u)^2)'.
\end{equation}
By \eqref{eq:As.10}, we get
$$Q_\lambda(u)-1=F_\lambda'(\cos u)^2+2F_\lambda'(\cos u)+\frac{\lambda}{4}\cdot(1+F_\lambda'(\cos u))\cos u\cdot\sin^2 2u.$$
 By the triangle inequality and  \eqref{eq:As.9n},
\begin{equation}\label{eq:As.20}
|Q_\lambda(u)-1|\leq 3|\lambda|(1+|\lambda|).
\end{equation}
Moreover, we get
\begin{equation*}
\begin{split}
Q_\lambda'(u)=
& 2(1+F_\lambda'(\cos u))\cdot F_\lambda''(\cos u)\cdot(- \sin u)\\
&+\frac{\lambda}{4}\cdot (1+F_\lambda'(\cos u))\cdot(-\sin u\sin ^22u+4\cos u\cdot\sin2u\cdot\cos 2u)\\
&
+\frac{\lambda}{4}\cdot \cos u\cdot\sin^22u\cdot F_\lambda''(\cos u)\cdot(-\sin u).
\end{split}
\end{equation*}
Since $F_\lambda''(\cos u)\cdot(-\sin u)=\frac{\lambda}{2}\cos u\cdot\sin 2u,$
it follows that
\begin{equation}\label{eq:As.21}
\begin{split}
Q_\lambda'(u)/\sin 2u=
& \lambda\cdot(1+F_\lambda'(\cos u))\left (\cos u-\frac{\sin u\sin 2u}{4}+\cos u\cos 2u\right)
\\
&+\frac{\lambda^2\cos^2 u\sin^22u}{8}
\end{split}
\end{equation}
for all $u\in(0,\pi)\setminus\{\pi/2\}.$
By applying the triangle inequality and \eqref{eq:As.9n} to \eqref{eq:As.21},
we have
\begin{equation}\label{eq:As.22}
|Q_\lambda'(u)/\sin 2u|\leq 4|\lambda|(1+|\lambda|)
\end{equation}
for all $u\in(0,\pi)\setminus\{\pi/2\}.$

By \eqref{eq:As.4}, \eqref{eq:As.9} and \eqref{eq:As.9n},
we obtain
$$x_\lambda'(u)^2+y'(u)^2=\sin^2u\cdot(1+F_\lambda'(\cos u) )^2+\alpha^2\cos^2u,$$
\begin{equation}\label{eq:As.23}
x_\lambda'(u)^2+y'(u)^2\leq(1+|F_\lambda'(\cos u)|)^2\leq(1+4|\lambda|/3)^2\leq 2(1+|\lambda|)^2
\end{equation}
and
\begin{equation*}
\begin{split}
(x_\lambda'(u)^2+y'(u)^2)^{'}
&=
\sin 2u \left(  (1+F_\lambda'(\cos u))^2+\lambda\sin^2u\cdot \cos u\cdot(1+F_\lambda'(\cos u))-\alpha^2\right)\\
&=\sin 2u\left(  1-\alpha^2+F_\lambda'(\cos u)^2+2F_\lambda'(\cos u)+\lambda\sin^2u\cos u(1+F_\lambda'(\cos u))\right).
\end{split}
\end{equation*}
Thus, by \eqref{eq:As.9n}
we obtain
\begin{equation}\label{eq:As.24}
|(x_\lambda'(u)^2+y'(u)^2)^{'}/\sin 2u-(1-\alpha^2)|\leq 4|\lambda|(1+|\lambda|)
\end{equation}
for all $u\in(0,\pi)\setminus\{\pi/2\}.$

From \eqref{eq:As.20}, \eqref{eq:As.22}, \eqref{eq:As.23} and  \eqref{eq:As.24},   it follows that
\begin{equation}\label{eq:As.27}
\left|\frac{Q'_\lambda(u)}{\sin 2u}\cdot (x_\lambda'(u)^2+y'(u)^2)\right|\leq 8|\lambda|(1+|\lambda|)^3,
\end{equation}
\begin{equation}\label{eq:As.26}
2\left|(Q_\lambda(u)-1)\cdot \frac{(x_\lambda'(u)^2+y'(u)^2)^{'}}{\sin 2u}\right|\leq 24|\lambda|(1+|\lambda|)^3,
\end{equation}
and 
\begin{equation}
\label{eq:As.25}
2(x_\lambda'(u)^2+y'(u)^2)^{'}/\sin 2u\geq 2(1-\alpha^2)-8|\lambda|(1+|\lambda|)
\end{equation}
for all $u\in(0,\pi)\setminus\{\pi/2\}.$

Then, we get the following lemma:
\begin{lemma}\label{lemAs.7}
There exists a number $\lambda_0\in(0,1)$ such that for each $\lambda\in(-\lambda_0,\lambda_0),$
the function $-m''_\lambda(x)/m_\lambda(x)$ is decreasing on $(0,r_e)$ and increasing on $(r_e,\pi).$
\end{lemma}
\begin{proof} Let us 
choose a number $\lambda_0\in(0,1)$ so as to satisfy
$$8|\lambda|(1+|\lambda|)^3+24|\lambda|(1+|\lambda|)^3+8|\lambda|(1+|\lambda|)    
<2(1-\alpha^2)$$
for all $\lambda$ with $|\lambda|<\lambda_0.$
We will check that the number $\lambda_0$ has the required property.
Choose any $\lambda\in(-\lambda_0,\lambda_0)$ and fix it.
From \eqref{eq:As.19}, it is clear that
$$P_\lambda(u)=Q'_\lambda(u)(x_\lambda'(u)^2+y'(u)^2)-2(Q_\lambda(u)-1)\cdot
(x_\lambda'(u)^2+y'(u)^2)'-2(x_\lambda'(u)^2+y'(u)^2)',$$
and
\begin{align*}
P_\lambda(u)/\sin2u &\leq 
\left|\frac{Q'_\lambda(u)}{\sin 2u}\cdot (x_\lambda'(u)^2+y'(u)^2)\right|
+2\left|(Q_\lambda(u)-1)\cdot \frac{(x_\lambda'(u)^2+y'(u)^2)^{'}}{\sin 2u}\right|\\
&
-\frac{2(x_\lambda'(u)^2+y'(u)^2)^{'}}{\sin 2u}.\\
\end{align*}
From \eqref{eq:As.27}, \eqref{eq:As.26}, and 
\eqref{eq:As.25}, it follows that $P_\lambda(u)/\sin 2u<0$ for all $u\in(0,\pi)\setminus\{\pi/2\}.$
Hence, $P_\lambda(u)<0$ on $(0,\pi/2)$ and $P_\lambda(u)>0$ on $(\pi/2,\pi).$ 
This implies,
due to \eqref{eq:As.18} and \eqref{eq:As.6}, that the function $-m_\lambda''(r)/m_\lambda(r)$ is strictly decreasing on $(0,r_e)$ and strictly increasing on $(r_e,\pi).$
$\qedd$
\end{proof}


\bigskip
\noindent

\medskip
\noindent
\begin{tabbing}
\hspace{100mm}\=\hspace{100mm}\kill
Minoru Tanaka\> Toyohiro Akamatsu\\
School of Science\>School of Science\\
Department of Mathematics\>Department of Mathematics\\
Tokai University, Hiratsuka City\>Tokai University, Hiratsuka City\\
Kanagawa Pref., 259\,--\,1292, Japan\>Kanagawa Pref., 259\,--\,1292, Japan\\
{\tt tanaka@tokai-u.jp}\>{\tt akamatsu@tokai-u.jp }
\end{tabbing}

\begin{tabbing}
\hspace{100mm}\=\hspace{100mm}\kill
Robert Sinclair\> Masaru Yamaguchi\\
Faculty of Economics\>School of Science\\
Hosei University\>Department of Mathematics\\

4342 Aihara-machi, Machida       
\>Tokai University, Hiratsuka City\\
Tokyo, 194\,--\,0298, Japan\>
Kanagawa Pref., 259\,--\,1292, Japan\\

{\tt sinclair.robert.28@hosei.ac.jp} \>{\tt ym24896@tsc.u-tokai.ac.jp}\\

\end{tabbing}

\end{document}